\documentclass[10pt,a4paper]{amsart}
\usepackage{t1enc}

\usepackage[utf8]{inputenc}
\usepackage[english]{babel}
\usepackage{graphicx}
\usepackage{amsmath}
\usepackage{amssymb}
\usepackage{amsthm}
\usepackage{mathrsfs}
\usepackage[all]{xy}
\usepackage{enumerate}

\theoremstyle{plain}
\newtheorem{thm}{Theorem}[section]
\newtheorem{prop}[thm]{Proposition}

\newtheorem{cor}[thm]{Corollary}
\newtheorem{lem}[thm]{Lemma}

\theoremstyle{definition}
\newtheorem{dfn}{Definition}[section]

\theoremstyle{remark}
\newtheorem*{rem}{Remark}

\DeclareMathOperator{\Aut}{Aut}

\DeclareMathOperator{\Gal}{Gal}

\DeclareMathOperator{\Spec}{Spec}

\DeclareMathOperator{\ga}{\pi_1^{alg}}
\DeclareMathOperator{\gtop}{\pi_1^{top}}

\DeclareMathOperator{\gtemp}{\pi_1^{temp}}

\DeclareMathOperator{\Covalg}{Cov^{alg}}
\DeclareMathOperator{\Covtop}{Cov^{top}}
\DeclareMathOperator{\Covtemp}{Cov^{temp}}

\DeclareMathOperator{\Ker}{Ker}

\DeclareMathOperator{\Hom}{Hom}

\DeclareMathOperator{\Set}{Set}

\DeclareMathOperator{\Stab}{Stab}

\DeclareMathOperator{\an}{an}

\DeclareMathOperator{\Gm}{\mbf G_m}

\DeclareMathOperator{\Frac}{Frac}
\DeclareMathOperator{\ord}{ord}
\DeclareMathOperator{\PGL}{PGL}
\DeclareMathOperator{\can}{can}

\DeclareMathOperator{\degtr}{degtr}
\DeclareMathOperator{\val}{val}

\DeclareMathOperator{\Conv}{Conv}

\newcommand{\ie}{\emph{i.e.} }
\newcommand{\findem}{\end{proof}}
\newcommand{\dem}{\begin{proof}}

\newcommand{\da}{\begin{displaystyle}}
\newcommand{\db}{\end{displaystyle}}

\newcommand{\mcal}{\mathcal}
\newcommand{\mbf}{\mathbf}
\newcommand{\bb}{\mathbb}
\newcommand{\mbb}{\mathbb}

\newcommand{\C}{\mathbf C}

\newcommand{\Q}{\mathbf Q}
\newcommand{\Z}{\mathbf Z}
\newcommand{\F}{\mathbf F}
\newcommand{\N}{\mathbf N}
\renewcommand{\P}{\mathbf P}
\newcommand{\A}{\mathbf A}
\newcommand{\G}{\mathbb G}
\newcommand{\T}{\mathbb T}

\title{Resolution of non-singularities for Mumford curves}
\author{Emmanuel Lepage}

\begin{document}
\maketitle 


\section*{Introduction}

In this article, given a hyperbolic curve $X$ over $\overline{\Q}_p$, we
are interested in finding a finite étale cover $Y$ of this curve such that
the stable reduction of the cover has irreducible components lying over the
smooth locus of the stable reduction of $X$. Such techniques of resolution
of nonsingularities are often used in anabelian geometry. We will apply our
results to the anabelian study of the tempered fundamental group.

In \cite[th. 0.2]{resnonsing}, A. Tamagawa proved that for every hyperbolic
curve $X=\overline X\backslash D$ and every closed point $x$ of the stable
reduction of $X$, the exists a finite étale cover $Y$ and an irreducible
component $y$ of the stable reduction of $Y$ lying above $x$.
We would like to generalize this to all the semistable reductions of $X$:
given a semistable model $\mathcal X$ of $X$ and a closed point $x$ of the
special fiber $\mathcal X_s$ of $X$, is there a finite étale cover $Y$ and an
irreducible component $y$ of the minimal semistable model $\mcal Y$ of $Y$
above $\mathcal X$ such that $y$ lies above $x$? To give an example of anabelian motivation for this kind of resolution of nonsingularities, as shown by F. Pop and J. Stix in \cite[cor. 41]{popstix}, if $X_0$ is a geometrically connected hyperbolic curve over a finite extension $K$ of ${\Q}_p$ such that $X_{0,\overline{\Q}_p}$ satisfies this kind of resolution of nonsingularities, every section of $\ga(X_0)\to \Gal(\overline{\Q}_p/K)$ has its image in a decomposition group of a unique valuative point. In this article, we will
prove that this resolution of nonsingularities is satisfied by Mumford curves.

First, let us translate this in an analytic setting. Let $X^{\an}$ be the
Berkovich space of $X$. Given a semistable model $\mcal X$ of $X$, there is
a reduction map $\pi_{\mcal X}:X^{\an}\to\mcal X_s$. If $\eta_y$ is the
generic point of an irreducible component of $\mcal Y$, the subset
$\pi_{\mcal Y}^{-1}(\eta_y)\subset Y^{\an}$ is reduced to a single
point. We denote by $V(\mcal Y)$ the set of points of $Y^{\an}$ whose image
by $\pi_{\mcal Y}$ is a generic point. Therefore our question reduces to
the following: is there $Y$ and an element of $V(\mcal Y)$ which is mapped
to $\pi_{\mcal X}^{-1}(x)$? Since $V(\mcal Y)$ contains $V(Y):=V(\mcal
Y_0)$ where $\mcal Y_0$ is the stable model of $Y$ and $\pi_{\mcal
  X}^{-1}(x)$ is a non-empty open subset of $X^{\an}$, it is enough to show
that the union $\widetilde V(X)$ of the images of $V(Y)$ in $X^{\an}$, where $Y$
runs over finite étale covers of $X$, is dense in $X^{\an}$.
\begin{thm}[{th. \ref{rnsmumford}}]
 Let $X$ be a hyperbolic Mumford curve. Then $\widetilde V(X)$ is dense in $X^{\an}$.
\end{thm}

To do so, we will study $\mu_{p^n}$-torsors of $X$. Projective systems of
$\mu_{p^n}$-torsors of $X$ are classified by $H^1(X,\Z_p(1))$. 
Let $c$ be an element of $H^1(X,\Z_p(1))$. Let $x$ be a $\mbf C_p$-point of
$X$. There is a small rigid neighborhood $D$ of $x$ in $X_{\C_p}$
isomorphic to a disk and a morphism $f:D\to\mbf G_{m}$, \ie an element
$f\in O^*(D)$, such that the restriction of $c$ to $D$ is the pullback of
the canonical element of $H^1(\G_m,\Z_p(1))$.  Let $Y_n\to X$ be the
$\mu_{p^n}$-torsor induced by $c$. For $n$ big enough, there is a smallest
closed disk $D_n$ of $D$ centered at $x$ such that $Y_{n|D_n}\to X_{n|D_n}$
is a non-trivial cover. Then the behavior of the restriction of $Y_n$ to
the Berkovich generic point $x_n$ of $D_n$ for $n$ big enough only depends
on the ramification index $e$ of $f:D\to \mbf G_m$ at $x$. More precisely,
if $y_n$ is a preimage of $x_n$ in $Y_n$, the extension $\mcal H(y_n)/\mcal
H(x_n)$ of complete residue fields induce an extension $k(y_n)/k(x_n)$ of
their reduction in characteristic $p$. This extension is an  Artin-Schreier
extension: $k(x_n)$ is isomorphic to $\overline{\F}_p(X)$ and $k(y_n)$ is
isomorphic to $k(x_n)[Y]/(Y^p-Y-X^e)$. The general study of Artin-Schreier
extensions tells us that if $e$ is not a power of $p$, $k(y_n)$ is not a
rational extension of $\overline{\F}_p$. This implies that $y_n\subset
V(Y_n)$ and that $x_n\in \widetilde V(X)$ and therefore $x$ lies in the closure of
$\widetilde V(X)$.

By Hodge-Tate theory, one has a canonical decomposition
$H^1(X,\Z_p(1))\otimes_{\Z_p} \C_p=H^1(X_{\C_p},O_{X_{\C_p}})(1)\oplus
H^0(X_{\C_p},\Omega_{X_{\C_p}})$. Let us consider the induced map
$p:H^1(X,\Z_p(1))\to H^0(X_{\C_p},\Omega_{X_{\C_p}})$.  
Assume now $X$ is a Mumford curve over $\overline{\Q}_p$. Then the image of $p$ lies
in $H^0(X,\Omega_X)$ and, for $c\in
H^1(X,\Z_p(1))$, the restriction of $p(c)$ to $D$ is $\frac{df}{f}$. Let $\Omega\subset \P^1$ be the universal topological
cover of $X$. If $x\in \Omega(\overline{\Q}_p)$, one can find a rational
function $f$ with no poles nor zero in $\Omega$ such that $\frac{df}{f}$
has a zero at $x$ with multiplicity $m$ such that $m+1$ is not a power of
$p$. Let $c_f$ be the pullback of the canonical element of
$H^1(\Gm,\Z_p(1))$ along $f:\Omega\to \Gm$, and let $x_n\in \Omega$ be
defined as previously.
Then for the topology of the uniform convergence on every compact on
$O^*(D)$, we will approximate $f$ by elements of $\bigcup_{X'} \Theta(X')$,
where $X'$ runs over finite topological pointed covers of $X$ and
$\Theta(X')$ is the set of theta functions of $X'$. Using this, one can
construct for every $n$ a finite topological cover $X'$ of $X$ and a
$\mu_{p^n}$-torsor $Y\to X'$ such that the preimage in $Y$ of the image in
$X'$ of $x_n$ lies in $V(Y)$. Therefore $x\in \widetilde V(X)$. Since
$X_{\overline{\Q}_p}$ is dense in $X^{\an}$, one gets the density of
$\widetilde V(X)$.

In a second part, we use the resolution of nonsingularities to study the
tempered fundamental group. One shows the following:
\begin{thm}[{th. \ref{resnonsing}}]\label{resnonsingintro}
Let $X_1$ and $X_2$ be two Mumford curves over $\overline{\Q}_p$. Given an
isomorphism between their tempered fundamental groups, there is a canonical
homeomorphism between there Berkovich spaces. 
\end{thm}
The strategy is the following. For every semistable model of a curve $X$,
there is a retraction from $\overline X^{\an}$ to the graph of this
semistable reduction. One gets a map from $\overline X^{\an}$ to the
projective limit of graphs of semistable reductions, which is a
homeomorphism. If $X$ is a Mumford curve, by resolution of
nonsingularities, semistable models of the form $\mcal Y/G$, where $Y$ runs
over finite Galois cover of $X$, $\mcal Y$ is the stable model of $Y$ and
$G=\Gal(Y/X)$, are cofinal among semistable models of $X$. However, a
theorem of S. Mochizuki tells us that one can recover the graph of the
stable reduction from the tempered fundamental group
(\cite[cor. 3.11]{mochi}). Therefore if $Y_1$ is a finite Galois cover of
$X_1$ and $Y_2$ is the corresponding finite Galois cover given by the
isomorphism of tempered fundamental groups, the graph $\mbb G_{Y_1}$ of
$\mcal Y_1$ is canonically isomorphic to the graph $\mbb G_{Y_2}$ of $\mcal
Y_2$, and one gets a similar isomorphism after quotienting by
$\Gal(Y_1/X_1)\simeq\Gal(Y_2/X_2)$. The problem is to recover from the tempered fundamental group the transition maps between
the geometric realisation of the different graphs.

At the end of the article, we will be interested by the anabelianness of the tempered fundamental group for punctured Tate curves:
\begin{thm}[{th. \ref{tatecurves}}]\label{tatecurveintro}
 Let $q_1,q_2\in\overline{\Q}_p$ such that $|q_1|,|q_2|<1$. Assume there
 exists an isomorphism $\psi$ between the tempered fundamental groups of
 $(\Gm/q_1^{\Z})\backslash\{1\}$ and $(\Gm/q_2^{\Z})\backslash\{1\}$. Then
 there exists $\sigma\in \Gal(\overline{\Q_p}/\Q_p)$ such that
 $q_2=\sigma(q_1)$.
\end{thm}
However, the proof of this result does not build any particular element of $\Gal(\overline{\Q_p}/\Q_p)$.

Let $\Omega_i=\Gm\backslash \{q_i^n\}_{n\in\Z}$ and $X_i=\Omega_i/q_i^{\Z}$. According to theorem \ref{resnonsingintro}, the isomorphism of tempered fundamental groups induces a homeomorphism $\bar\psi:\Gm^{\an}\to \Gm^{\an}$ which maps $q_1^n$ to $q_2^n$ for every $n\in \Z$. 

Elements of $O(\Omega_i)^*$ correspond, up to a scalar, to a current on the semitree $\mbb T_i$ of $\Omega_i$. Since  $\psi$ induces an isomorphism $\mbb T_1\simeq \mbb T_2$, one gets a group isomorphism $\alpha:O(\Omega_1)^*\to O(\omega_2)^*$. The crucial point will be to show that for every $f\in O^*(\Omega_1)$ and $z\in \Gm(\C)$, the multiplicity of $\frac{df}{f}$ at $z$ coincide with the multiplicity of $\frac{d\alpha(f)}{\alpha(f)}$ at $\bar\psi(z)$. By density of $\Z$ in $\Z_p$, one also gets a similar result for $\Z_p$-linear combinations of differential $1$-forms as $\frac{df}{f}$. 
Once one knows this,  one can build, for every $n\in\N$, an element $f\in O^*(\Omega_1)$ such that $\frac{df}{f}(1)=q_1^n$ and $\frac{d\alpha(f)}{\alpha(f)}(1)=q_2^n$: one therefore gets that for every polynomial $P\in \Z_p[X]$, $P(q_1)=0$ if and only if $P(q_2)=0$.

\section{Berkovich geometry of curves}

In the following $K=\mbf C_p$, $k$ is its residue field, which is
isomorphic to $\overline{\mbf F}_p$. The norm will be chosen so that
$|p|=p^{-1}$ and the valuation so that $v(p)=1$. 
All valued fields will have valuations with values in $\mathbf{R}_{\geq 0}$.

If $X$ is an algebraic variety over $K$,
one can associate to $X$ a topological set $X^{\an}$ with a continuous map
$\phi:X^{\an}\to X$ defined in the following way.
A point of $X^{\an}$ is an equivalence class of morphisms $\Spec K'\to X$
over $\Spec K$ where
$K'$ is a complete valued extension of $K$. Two morphisms
$\Spec K'\to X$ and $\Spec K''\to X$ are equivalent if there exists a common
valued extension $L$ of $K'$ and $K''$ such that
\[
\xymatrix{\Spec L \ar[r] \ar[d] & \Spec K''\ar[d] \\ \Spec K' \ar[r] &
  X}
\]
commutes. In fact, for any point $x\in X^{\an}$, there is a unique smallest
such complete valued field defining $x$ denoted by $\mcal H(x)$ and called the completed residue field of $x$.
We denote by $k(x)$ the residue field of $\mcal H(x)$ and by $\val(x)\subset \mbf R_{>0}$ the group of values of $\mcal H(x)$. Forgetting the valuation, one gets points $\Spec(K)\to X$ from the same equivalence class of points: this defines a point of $X$, hence the map $X^{\an}\to X$.
If $U=\Spec A$ is an affine open subset of $X$, every $x\in\phi^{-1}(U)$
defines a seminorm $|\ |_{x}$ on $A$. The topology on $\phi^{-1}(U)$ is
defined to be the coarsest such that $x\mapsto |f|_x$ is continuous for
every $f\in A$.

The space $X^{\an}$ is locally compact, and even compact if $X$ is proper.
In fact $X^{\an}$ is more than just a topological space: it can be enriched
into a $K$-analytic
space, as defined by Berkovich in~\cite{berk}.

\smallskip

Points of $\A^{1,an}$ are of four different types and are described in the following way:
\begin{itemize}
\item A closed ball $B=B(a,r) \subset \C_p$ of center $a$ and radius $r$ defines a point $b=b_{a,r}$ of $\A^{1,\an}$ by 
\[
|f|_b=\sup_{x\in B} |f(x)|=\max_{i\in \N} |a_i|r^i\text{ if } f=\sum_{i\in \N}a_iX^i.
\]
The point $b_{a,r}$ is said to be of type 1 if $r=0$, of type 2 if $r\in p^{\Q}$ and of type 3 otherwise. The pairs $(a,r)$ and $(a',r')$ define the same point if and only if $B(a,r)=B(a',r')$, \ie $r=r'$ and $|a-a'|\leq r$.

\item A decreasing family of balls $E=(B_i)_{i\in I}$ with empty intersection
   defines a point by 
   \[
   |f|_E=\inf_{i\in I} |f|_{b_i}.
   \]
    Such a point is said to be of type 4.
\end{itemize}
If $r\in p^{\Q}$ and $a\in \C_p$ are such that $|a|=r$, then $|X/a|_{b_{a,r}}=1$ and $k(b_{a,r})=k(\overline{X/a})$.

The classification by type of points can be generalized to curves in such a way that it is preserved by finite morphisms: a point $x$ is of type:
\begin{itemize}
 \item 1 if $\mcal H(x)=\C_p$;
\item 2 if $\degtr k(x)/\overline{\F}_p=1$;
\item 3 if $\val(x)\neq p^{\mbf Q}$;
\item 4 otherwise (\ie $\mcal H(x)/\C_p$ is an immediate extension).
\end{itemize}
If $x$ is of type $2$, we denote by $g_{k(x)}$ the genus of the proper $\overline{\F}_p$-curve whose field of fraction is $k(x)$.

\smallskip

Let $X$ be a proper and smooth $K$-curve.
The topological space is a quasipolyhedron in the sense of \cite[def. 4.1.1]{berk}:
there exists a base of open subsets $U$ such that:
\begin{itemize}
 \item $\overline U\backslash U$ is finite;
 \item $U$ is countable at infinity;
 \item for every $x\neq y\in U$, there exists a unique subset $[x,y]\subset U$ such that homeomorphic to $[0,1]$ with endpoints $x$ and $y$.
\end{itemize}
A quasipolyhedron that satisfies itself these three porperties is said to be simply-connected.

If $X$ is a curve, the set of points of type different from $2$ is totally disconnected. Therefore, on every nonconstant path, there exists infinitely many points of type $2$.

The topological universal cover $X^{\infty}$ is a simply-connected quasipolyhedron.
Therefore any subset $I$ of $X^\infty$ is contained in a smallest connected
subset $\Conv(I)=\bigcup_{(x,y)\in I^2}[x,y]$. If $I,J$ are closed connected subsets
of $X^{\infty}$, one denotes by $[I,J]=\bigcap_{(x,y)\in I\times
  J}[x,y]$. It is homeomorphic to $[0,1]$ if $I\cap J=\emptyset$ and $I\cap
[I,J]$ and $J\cap [I,J]$ are reduced to a point. If $x\in X^{\infty}$ and
$I$ is a closed subset of $X^\infty$, one denotes by $r_I(x)$ the unique
element of $[x,I]\cap I$. The map $r_I:X^{\infty}\to I$ is a continuous
retraction of the embedding $I\to X^{\infty}$.

Let $\mcal X$ be a semistable model.
There is a specialization map $\pi_{\mcal X}:X^{\an}\to \mcal X_k$ defined in the
following way. If $x\in X^{\an}$, the morphism $\Spec \mcal H(x)\to X$
extends in a morphism $\Spec O_{\mcal H(x)}\to\mcal X$ by properness of
$\mcal X\to \Spec O_K$, hence a morphism $\Spec k(x)\to \mcal X_k$: the
image of this morphism is $\pi_{\mcal X}(x)$. This specialization map is anticontinuous: the preimage of a closed subset is an open subset.

If $z$ is the generic point of an irreducible component of $\mcal X_k$,
then $z$ is of codimension $1$ in $\mcal X$ and thus $O_{\mcal X,z}$ is a
valuation ring. The completion of $\Frac O_{\mcal X,z}$ defines a point
$b_z$ of $X^{\an}$ which is the unique element of $\pi_{\mcal X}^{-1}(z)$. One
denotes by $V(\mcal X)\subset X^{\an}$ the set of such $b_z$ and by
$V(\mcal X)^{\infty}$ its preimage in $X^{\infty}$. 

One has $X^{\an}\backslash V(\mcal X)=\coprod \pi_{\mcal X}^{-1}(x)$ where $x$ goes through closed points of $\mcal X_k$, and $\pi_{\mcal X}^{-1}(x)$ is open by anticontinuity of $\pi_{\mcal X}$. 
In particular, if $z,z'\in X^{\an}$ are such that $\pi_{\mcal X}(z)\neq\pi_{\mcal X}(z')$, then every path joining $z$ to $z'$ meets $V(\mcal X)$.

One denotes  by $S(\mcal X)^{\infty}=\Conv(V(\mcal X)^{\infty})$ and by $r^\infty_{\mcal X}$ the
retraction $r_{S(\mcal X)^\infty}:X^\infty\to S(\mcal X)^\infty$ of the
embedding $\iota_{\mcal X}^\infty:S(\mcal X)^\infty\to X^\infty$. Since
$V(\mcal X)^\infty$ is $\Gal(X^{\infty}/X)$-invariant, $S(\mcal X)^\infty$
is also $\Gal(X^{\infty}/X)$-invariant and $r^\infty_{\mcal X}$  is
$\Gal(X^\infty/X)$-equivariant. One denotes by $S(\mcal X)$ the image of
$S(\mcal X)^{\infty}$ in $X^{\an}$: it is called the skeleton of $\mcal
X$. One denotes by $r_{\mcal X}$ the retraction $X^{\an}\to S(\mcal X)$
induced by $r^\infty_{\mcal X}$. The space $S(\mcal X)$ is compact and the
inclusion map $\iota_{\mcal X}:S(\mcal X)\to X^{\an}$ is a homotopy
equivalence. In fact $\mcal X$ is characterized by the fact that it is the
smallest subset $S$ of $X^{\an}$ that contains $V(\mcal X)$ and such that
$S\to X^{\an}$ is a homotopy equivalence.

If $z$ is a node of $\mcal X_k$, then $\pi^{-1}(z)$ is an open annulus. It contains a unique closed connected subset $S_z$ homeomorphic to $\mbf R$. More precisely, if one choses an isomorphism of analytic spaces $\pi^{-1}(z)\simeq \{z\in \A^{1,\an}|r_0<|T|_z<1\}$, then $S_z=\{b_{0,r}, r_0<r<1\}$ (in particular the points of type $2$ of $S_z$ can be identified with $\Q\cap (r_0,1)$).
If $z$ is a closed point of the smooth locus of $\mcal X_k$, then $\pi^{-1}(z)$ is an open disk.
Since, for every point $b$ of type $2$ of disks and annuli, $k(b)$ is a rational extension of $k$, one gets that if $x\in X^{\an}$ is a point of type $2$ such that $g_{k(x)}\neq 0$, then $x\in V(\mcal X)$.

One recovers $S(\mcal X)$ as the union of $V(\mcal X)$ and of $S_z$ for every node $z$ of $\mcal X_k$. One gets that $S(\mcal X)$ is homeomorphic to the dual graph of $\mcal X_k$.

If $X_1\to X_2$ is a morphism of proper and smooth $K$-curves and $\mcal X_1\to\mcal X_2$ is an extension to semistable $O_K$-models, then
\[\xymatrix{X_1^{\an}\ar[r]\ar[d] & \mcal X_{1,k}\ar[d]\\
X_2^{\an}\ar[r]  &\mcal X_{2,k}}
\]
is commutative.

Assume now that $\mcal X_1\to \mcal X_2$ is a morphism of semistable models of a same curve $X$.
If $z_1$ is the generic point of an irreducible component of $\mcal
X_{1,k}$ which maps to the generic point $z_2$ of an irreducible component
of $\mcal X_{2,k}$, the previous diagram tells us that
$b_{z_1}=b_{z_2}$. Since $\mcal X_1\to\mcal X_2$ is surjective, one gets
that $V(\mcal X_2)\subset V(\mcal X_1)\subset S(\mcal X_1)$. Since $S(\mcal
X_2)$ is the smallest subset of $X^{\an}$ that contains $V(\mcal X_2)$ and
such that $S(\mcal X_2)\to X^{\an}$ is a homotopy equivalence, one gets
that $S(\mcal X_2)\subset S(\mcal X_1)$. Similarly, one has $V(\mcal
X_2)^\infty\subset V(\mcal X_1)^\infty$ and $S(\mcal X_2)^\infty\subset
S(\mcal X_1)^\infty$.

Therefore, for every $x\in X^{\an}$, $[x,S(\mcal X_1)]\subset [x,S(\mcal X_2)]$, and thus $\iota_{\mcal X_1}r_{\mcal X_1}\in [x,S(\mcal X_2)]$. This implies that $r_{\mcal X_2}\iota_{\mcal X_1}r_{\mcal X_1}=r_{\mcal X_2}$.
Therefore, the maps $r_{\mcal X_1/\mcal X_2}:=r_{\mcal X_2}\iota_{\mcal X_1}:S(\mcal X_1)\to S(\mcal X_2)$ are compatible with composition and $(S(\mcal X))_{\mcal X}$ is a projective system of topological spaces. The maps $(r_{\mcal X})$ induce a conitnuous map \[r_X:X^{\an}\to \varprojlim_{\mcal X} S(\mcal X).\]
\begin{prop}\label{homeosst}
The map $r_X$ is a homeomorphism.
\end{prop}
\dem
Since $X^{\an}$ is compact, the surjectivity of $r_X$ follows from the surjectivity of each $r_{\mcal X}$.

Let $x\neq x'$ be such that $r_{\mcal X_0}(x)=r_{\mcal X_0}(x')$ where $\mcal X_0$ is the minimal model of $X$. Let $U$ be a simply connected open neighborhood of $r_{\mcal X_0}(x)$ in $S(\mcal X_0)$ and $V=r_{\mcal X_0}^{-1}(U)$. Since $V$ is a simply connected quasipolyhedron, there is a minimal connected subset $[x,x']\subset V$ containing $x$ and $x'$. It is homeomorphic to $[0,1]$ and has a natural order that make $x$ the smallest element. Let $x_1<x_2\in [x,x']$ be points of type $2$.

Since $x_i$ is of type $2$, $V\backslash \{x_i\}$ has infinitely many components. Since $\overline{\Q}_p$-points are dense in $X^{\an}$ one can find $z_{i,1},z_{i,2},z_{i,3}\in V\cap X(\overline{\Q}_p)$ lying in different connected components of $X^{\an}\backslash \{x_i\}$. Let $\mcal X$ be the stable model of the marked curve $(X,\{z_{ij}\}_{i=1,2;j=1,2,3}$. Since $\pi_{\mcal X}(z_{i,1})\neq \pi_{\mcal X}(z_{i,2})$, one has $S(\mcal X)\cap [z_{i,1},z_{i,2}]\neq\emptyset$ and therefore $r_{\mcal X}(z_{i,1})\in [z_{i,1},z_{i,2}]$.

Therefore, replacing $z_{i,2}$ by $z_{y,3}$, one gets that $r_{\mcal X}(z_{i,1})\in [z_{i,1},z_{i,2}]\cap [z_{i,1},z_{i,3}]=[z_{i,1},x_i]$ and similarly $r_{\mcal X}(z_{i,2})\in [z_{i,2},x_i]$. Since $S(\mcal X)$ is connected and intersects $[z_{i,1},x_i]$ and $[z_{i,2},x_i]$, $x_i\in S(\mcal X)$. Therefore, in $[x,x']$, $r_{\mcal X}(x)\leq x_1<x_2\geq r_{\mcal X}(x')$, which proves the injectivity of $r_X$. 

Since $X^{\an}$ is compact and $\varprojlim_{\mcal X} S(\mcal X)$ is Hausdorff, $r_X$ is a homeomorphism.
\findem

Let $\mcal X_1\to\mcal X_2$ be a morphism of semistable models of $X$. Let $z$ be a node of $\mcal X_{2,k}$. We will the notation \begin{equation}\label{eq:Amodel} A_{z,\mcal X_1}:=V(\mcal X_1)\cap S_z.\end{equation}
If one chooses an orientation of $S_z\simeq \mbf R$, $A_{z,\mcal X_1}$ then becomes an ordered set.
One denotes by $A_{z}=\bigcup_{\mcal X_1}A_{z,\mcal X_1}\subset S_z$.

\section{Resolution of non-singularities}
\subsection{definition}
Let $X=\overline X\backslash D$ be a hyperbolic curve over $K$.
Let $X^{\an}_{(2)}\subset X^{\an}$ be the subset of type $(2)$ points.

Let $\widetilde V(X)$ be the set of points $x$ of $X^{\an}$ such that there
exists a finite étale cover $f:Y\to X$ and $y\in V(Y)$ such that
$f(y)=x$. If $x\in\widetilde V(X)$, then $Y$ can be chosen to be Galois, so
that in particular $f^{-1}(x)\subset V(Y)$, since $V(Y)$ is Galois
equivariant. Let $\overline V(X)$ be the closure of $\widetilde V(X)$ in
$\overline X$. If $f:Y\to X$ is a finite étale cover $\widetilde
V(Y)=f^{-1}(\widetilde V(X))$ and $\overline V(Y)=f^{-1}(\overline V(X))$.
One has $\widetilde V(X)\subset X^{\an}_{(2)}$.
\begin{dfn} One says that $X$ satisfies resolution of non-singularities, or $RNS(X)$ for short, if 
$\widetilde V(X)= X^{\an}_{(2)}$.
\end{dfn}
\begin{prop}\label{rnsprop}
Let $X=\overline X\backslash D$ be a curve. The following are equivalent:

\begin{enumerate}
 \item $\widetilde V(X)=X^{\an}_{(2)}$;
\item $\overline V(X)=\overline X^{\an}$;
\item $X(K)\subset\overline V(X)$.
\end{enumerate}
\end{prop}
\dem
\begin{itemize}
 \item $(i)\Rightarrow (ii)$. Points of type $(2)$ are dense in $X_{\an}$.
\item $(ii)\Rightarrow (iii)$ is obvious.
\item $(iii)\Rightarrow (i)$. Let $x\in X_{(2)}$. Then
  $X^{\an}\backslash\{x\}$ has infinitely many components and they are open
  in $X^{\an}$. Since $X(K)$ is dense in $X^{\an}$, each of this components
  intersect $\widetilde V(X)$. Let $x_1, x_2,x_3\in \widetilde V(X)$ lying
  in different connected components $S_1,S_2,S_3$ of
  $X^{\an}\backslash\{x\}$. Let $f:Y\to X$ be a finite cover such that
  there is $y_i$ over $x_i$ lying in $S(Y)$ for $i=1,2,3$. Up to
  replacing $Y$ by a Galois closure, one can assume that $Y\to X$ is
  Galois. Since the image $T$ of $S(Y)$ in $X^{\an}$ is connected and
  $x_1x_2\in T$, one has $x\in T$.  Since $S(Y)\subset Y^{\an}$ is
  $\Gal(Y/X)$-invariant, $S(Y)=f^{-1}(T)$. Let $y\in f^{-1}(x)$. For
  any neighborhood $U$ of $x$, $U\cap T\cap S_i\neq \emptyset$. Thus for
  any neighborhood $V$ of $y$ small enough (for example, such that
  $f^{-1}(x)\cap V=\{y\}$), then $V\cap S(Y)\backslash\{y\}$ has at
  least three connected components. Thus $y\in V(Y)$ and
  $x\in\widetilde V(X)$.
\end{itemize}
\findem
If $X$ is a curve over $\overline{\Q}_p$, then $X(\overline{\Q}_p)$ is dense in $X(\overline{\C}_p)$, thus $X$ satisfies resolution of nonsingularities if and only if $X(\overline{\Q}_p)\subset\overline V(X)$.

If $Y\to X$ is a morphism of hyperbolic curves over $\overline{\Q}_p$ and $\mcal X$ is a semistable model of $\overline X$, there exists a minimal semistable model $\mcal Y$ of $Y$ above $\mcal X$ ($\mcal Y$ is the stable marked hull of the normalization of $\mcal X$ in $K(Y)$, in the sense of \cite[cor. 2.20]{liucompositio}).

\begin{prop}
Let $X$ be a $\overline{\Q}_p$-curve which satisfies resolution of nonsingularities. Let $\mcal X$ be a semistable model of $\overline X$ and let $x$ be a closed point of $\mcal X_k$. There exists a finite cover $Y\to X$ such that $\mcal Y_k$ has a vertical component above $x$, where $\mcal Y$ is the minimal semistable model of $Y$ above $\mcal X$.
\end{prop}
\dem
Let $\pi_{\mcal X}:\overline X^{\an}\to\mcal X_k$ be the specialization map. Since $x$ is closed, $\pi_{\mcal X}^{-1}(x)$ is open, and therefore contains a point $\tilde x$ of type $2$. Let $Y$ be a cover of $X$ and let $\tilde y\in V(Y)$ be above $\tilde x$. Then, for the stable model $\mcal Y_0$ of $Y$, $\tilde y$ specializes via $\pi_{\mcal Y_0}$ to a generic point of $\mcal Y_{0,k}$. Therefore, $\tilde y$ specializes to a generic point for every semistable model of $Y$, in particular for $\mcal Y$. Then $y:=\pi_{\mcal Y}(\tilde y)$ is mapped to $x$ and therefore the closure of $y$ is a vertical component above $x$.
\findem

\subsection{Splitting points of $\mbf Z_p(1)$-torsors}
The map $\Gm\stackrel{(\ )^n}{\to}\Gm$ defines a $\mu_n$-torsor over $\Gm$. The corresponding element of $H^1(\Gm,\mu_n)$ is denoted by $c_{\can,n}$.
 
Let $D$ be a disk centered at $0$ and let $f:D\to\Gm$ be a non constant morphism. Let $c_n=f^*c_{\can,p^n}\in H^1(D,\mu_{p^n})$. Let $Y_n\to D$ be the corresponding $\mu_{p^n}$-torsor.

Let $f(X)=\sum_{k\geq 0} a_kX^k$ be the power series of $f$.
Let $e_0(f)=\inf \{k\geq 1|a_k\neq 0\}$ be the ramification index of $f$ at $0$.

Let $r_0(c_n)=\inf\{r>0|Y_n\text{ is not split above } b_{0,r}\}$ when it
exists. Then $Y_n$ is trivial above  $D(0,r_0(c_n)^-)$ (otherwise it could
be extended in a non trivial finite cover of $\P^1$).

Let $y_n\in Y_n$ be above $x_n=b_{0,r_0(c_n)}$. The cover $Y_n\to D$
induces a morphism $\mcal H(x_n)\to\mcal H(y_n)$ of complete valued
field. Let $k(x_n)\to k(y_n)$ be the morphism of their residue fields.

We want to study the asymptotic behavior of $r_0(c_n)$ when $n$ goes to
$\infty$ and $k(y_n)$.
\begin{prop}\label{rayonconv}Their exists $C$ such that, for $n$ big
  enough, $r_0(c_n)=Cp^{-\frac{n}{e_0(f)}}$.

Moreove, for $n$ big enough, $[k(y_n):k(x_n)]=[\mcal H(y_n):\mcal
H(x_n)]=p$ and $k(y_n)$ is isomorphic to $k(X)[T]/(T^p-T-X^{e_0(f)})$.
\end{prop}
\dem
Up to multiplying $f$ by a constant, one can assume $f(0)=1$. Let
$N=e_0(f)$. Up to replacing $D$ by a smaller disk, one can assume that
$\sum_{k\geq N} a_kX^k$ has a $N$th root $t$, so that $f=1+t^N$. Moreover
up to replacing $D$ by a smaller disk, one can assume that $t$ induces an
isomorphism $t:D\to D_0$ where $D_0$ is also a disk centered at $0$. Since
$t$ maps $b_{0,r}$ to $b_{0,\lambda r}$ where $\lambda$ is a constant, it
is enough to prove the result for $ft^{-1}$. One can thus assume that
$f=1+X^N$.

One has $f_{n-1}:=f^{1/p^{n-1}}=\sum_k b_kX^{Nk}$ where
$b_k=\binom{1/p^{n-1}}{k}$ and $v_p(b_k)=-k(n-1)-v_p(k!)$. The series
$f^{1/p^{n-1}}$ is convergent on the disk of radius
$p^{-(n-1+\frac{1}{p-1})/N}$. By replacing $n-1$ by $n$, one gets that
$r_0(c_n)\geq \lambda_n:=p^{-(n+\frac{1}{p-1})/N}$. Let $1+y\in O(Y_n)$ be
the $p^n$th root of $f$ such that $y(0)=0$.
Then $y$ satisfies the equation 
\begin{equation}\label{eq1} \sum^p_{k=1} \binom{p}{k}y^k=\sum_{k\geq 1} b_kX^{Nk}. \end{equation}
Let $b=b_{0,\lambda_n}\in D$. Let $b'$ be above $b$ in $Y_n$ and let $b''$
be the image of $b'$ in $Y_{n-1}$. Since $\lambda_{n-1}>\lambda_n$, the
torsor $c_{n-1}$ is split at $b''$, one has $[\mcal H(b'):\mcal
H(b)]=[\mcal H(b'):\mcal H(b'')]|p$ and $\mcal H(b)=\mcal
H(b')/((1+y)^p-f_{n-1})$.

\smallskip

At $b$, one has $|b_1X^N|_b=|\frac{1}{p^{n-1}}t^N|_b=p^{-\frac{p}{p-1}}$, and all the other terms in the right member of \eqref{eq1} have smaller norms: $|b_kt^{Nk}|_b=p^{-\frac{kp}{p-1}+v_p(k!)}< p^{-k}\leq p^{-\frac{p}{p-1}}$ for every $k\geq 2$. In particular
\[|\sum^p_{k=1} \binom{p}{k}y^k|_{b'}=|\sum_{k\geq 1} b_kX^{Nk}|_b=p^{-\frac{p}{p-1}}.\]

If $|y|_{b'}<p^{-\frac{1}{p-1}}$, then $|\binom{p}{k}y^k|_{b'}=p^{-1}|y^k|<p^{-\frac{p}{p-1}}$ if $1\leq k<p$ and $|\binom{p}{k}y^k|_{b'}=|y|^p<p^{-\frac{1}{p-1}}$ if $k=p$, which is impossible since $|\sum^p_{k=1} \binom{p}{k}y^k|_{b'}=p^{-\frac{p}{p-1}}$. If $|y|_{b'}>p^{-\frac{1}{p-1}}$, then $|\binom{p}{k}y^k|_{b'}<|y^p|$ for every $1\leq k<p$ and therefore $|\sum^p_{k=1} \binom{p}{k}y^k|_{b'}=|y^p|>p^{-\frac{p}{p-1}}$, which is impossible.
Therefore, $|y|_{b'}=p^{-\frac{1}{p-1}}$. One gets that $|py|=|y^p|=p^{-\frac{p}{p-1}}$ and $|\binom{p}{k}y^k|_{b'}<p^{-\frac{p}{p-1}}$ for every $2\leq k\leq p-1$.

Therefore, in the ring $\{a\in\mcal H(b')||a|\leq p^{-\frac{p}{p-1}}\}/\{a\in\mcal H(b')||a|< p^{-\frac{p}{p-1}}\}$, equation \eqref{eq1} becomes $py+y^p=X^N/{p^{n-1}}$.

Let $z=y/a_1\in \mcal H(b')$ and $u=X/a_2\in\mcal H(b)\subset\mcal H(b')$ where $a_1^{p-1}=-p$ and $a_2^N=p^{n-1}a_1^p$. One has $|z|_{b'}=|u|_{b'}=1$ and $k(b)=k(\bar u)$. Let $\bar z$ and $\bar u$ be the classes of $z$ and $u$ in $k(b')$. Equation \eqref{eq1} induces in $k(b')$ the equality:
\[\bar z^p-\bar z=\bar u^N.\]
Therefore $k(b)[\bar z]\subset k(b')$ is a non trivial extension of $k(b)$. Since $[\mcal H(b'):\mcal H(b)]|p$, one gets that $[\mcal H(b'):\mcal H(b)]=p$, $b=x_n$, and $k(b')=k(b)[\bar z]$ is the wanted Artin-Schreier extension of $k(b)$ (because $k(b)=k(\bar u)$).
\findem
If one writes $e_0(f)=p^md$ where $d$ is prime to $p$, the genus of the Artin-Schreier curve $T^p-T=X^{e_0(f)}$ is $g=(d-1)(p-1)/2$ (cf. \cite[\S 2.2, eq. (8)]{lauderAS}).
In particular, if $e_0(f)$ is not a power of $p$, then, for $n$ big enough, $k(y_n)$ is not isomorphic to $k(X)$, so that if $Y_n$ is an anlytic domain of a curve $\widetilde Y_n$, one must have $y_n\in V(\widetilde Y_n)$.
It should be noticed that $e_0(f)=\ord_0(\frac{df}{f})+1$, where $\ord_0\omega$ denotes the $x$-adic valuation of $\frac{df}{fdX}\in \C_p[[X]]$.

\subsection{Resolution of non-singularities for Mumford curves}
In this subsection we show that Mumford curves over $\overline{\Q}_p$ satisfy resolution of nonsingularities.

A proper curve $X$ over $\overline{\Q}_p$ is a \emph{Mumford curve} if the following equivalent properties
are satisfied:
\begin{itemize}
\item all 
normalized irreducible components of its
stable reduction are isomorphic to $\P^1$,
\item $X^{\an}$ is locally isomorphic to $\P^{1,\an}$,
\end{itemize}

The universal topological covering $\Omega$ of $X^{\an}$ for a Mumford curve $X$ is an open subset of $\P^{1,\an}$. More precisely there is a
Shottky subgroup $\Gamma$ of $\PGL_2(\C_p)$, \ie a free finitely generated
discrete subgroup of $\PGL_2(\C_p)$, such that $\Omega=\P^{1,\an}\backslash\mcal L$ where $\mcal L$ is the closure of the set of $\C_p$-points stabilized by some nontrivial element of $\Gamma$. The points 
of $\mcal L$ are of  type 1, \ie are $\C_p$-points.
 Then $X$ is $p$-adic analytically uniformized as 
 \[
 X^{\an}=\Omega/\Gamma
 \]
 and $\Gamma=\gtop(X)$.

\begin{lem}
Let $I$ be an infinite subset of $K_0$, where $K_0\subset K$ is a finite extension of $\mbf Q_p$. Let $x\in\overline{\mbf Q}_p$ be a point not belonging to $I$. Let $E=\{k\in \N|\exists a\in \mbf Q_p^{(I)},k=\ord_x(\sum_{i\in I} \frac{a_i}{X-i})\}$ and $u_n=\#(E\cap [0,n])$ . Then the sequence $(u_n/n)_n$ does not go to $0$ when $n$ goes to infinity.
\end{lem}\label{densityord}
\dem
Up to replacing $K_0$ by $K_0[x]$ and $i$ by $i+x$ for every $i\in I$, one can assume $x=0$.
Let $V$ be the $\mbf Q_p$-vector subspace of $K_0[[X]]$ generated by $(\frac{1}{X-i})_{i\in I}$.
Let $C=[K_0:\mbf Q_p]$.

Let $V_n$ be the image of the $\mbf Q_p$-linear map $\phi_n:V\to K_0[[X]]/X^n$.
One has $n\in E$ if and only if $\Ker(\phi_n)\subsetneq\Ker(\phi_{n-1})$. Therefore \begin{itemize}
\item $\dim_{\mbf Q_p}V_n=\dim_{\mbf Q_p}V_{n-1}$ if $n\notin E$;
\item $\dim_{\mbf Q_p}V_n\leq C+\dim_{\mbf Q_p}V_{n-1}$.
\end{itemize}
Therefore $\dim_{\mbf Q_p}V_n\leq Cu_n$. However the morphism $f_n:V_n\otimes_{\mbf Q_p} K\to K[[X]]/X^n$ is surjective. Indeed let $i_1,\dots ,i_n$ be $n$ different elements of $I$ and let $\overline P\in K[[X]]/X^n$. Let $P$ be a representative of $\overline P$ in $K[X]$. Let $R$ be the reminder of $P\prod_{k=1}^n(X-i_k)$ by the division by $X^n$. Then, since $\deg(R)<n$, $\frac{R}{\prod_k(X-i_k)}=\sum_k\frac{a_k}{X-i_k}$ for some $(a_k)_k$ in $K^n$, and, since $\prod_k(X-i_k)$ is invertible in $K[[X]]$, $\overline P=f_n(\sum_k\frac{1}{X-i_k}\otimes a_i)$, which proves that $f_n$ is surjective. Therefore $\dim_{\mbf Q_p} V_n=\dim_{\mbf K}V_n\otimes_{\mbf Q_p} K\geq n$ and $u_n/n\geq 1/C$. 
\findem
Let $I$ be an infinite subset of $K_0$, where $K_0\subset K$ is a finite extension of $\mbf Q_p$. Let $x\in \overline{\mbf Q}_p$.
Lemma \ref{densityord} shows that there exists $(a_i)\in \mbf Q_p^{(I)}$ be such that $g=\sum\frac{a_i}{X-i}$ is such that $\ord_x(g)+1$ is not a power of $p$. Up to multiplying all the $a_i$ by $p^{\max_{i\in I}(-v_p(a_i))}$, one can assume $a_i\in\mbf Z_p$. Let $I_0\subset I$ be the support of the family $(a_i)_{i\in I}$.

For $i\in I_0$ and $n\geq 0$, let $a_{i,n}\in \mbf Z$ be such that $v_p(a_{i,n}-a_i)\geq n$.

Let $f_n=\prod_{i\in I_0} (\frac{X-i}{x-i})^{a_{i,n}}:\mbf P^1\backslash I_0\to \Gm$. Let $D=\{z\in \Gm||X-x|_z<\min_{i\in I_0}|x-i|\}$. The sequence $(f_n)$ is uniformly convergent on every affinoid subset of $D$ and defines over $D$ a morphism  $f:D\to\mbf G_m$ and $f'/f=g$ over $D$.

Let $c_n: Y_n\to \mbf P^1\backslash I_0$ be the $\mu_{p^n}$-torsor over $\mbf P^1\backslash I_0$ obtained by pulling back along $f_n$ the canonical torsor (remark that $c_n$ only depends on $a_i$ and not on $a_{i,n}$).
The restriction of $c_n$ to $D$ is also the pullback of the canonical torsor along $f$.

According to proposition \ref{rayonconv}, there is a point $y_n\in Y_n$ such that $g_{k(y_n)}\geq 1$ and $c_n(y_n)\to x$.

One gets the following result:
\begin{prop}\label{torsensinf}
Let $I$ be an infinite subset of a finite extension of $\mbf Q_p$. For every $x\in \overline{\mbf Q}_p$, there is a finite subset $I_0\subset I\backslash\{x\}$ and a $\mbf Z_p(1)$-torsor $c=(c_n:Y_n\to \mbf P^1\backslash I_0)$ of $\mbf P^1\backslash I_0$ and for every $n\geq 1$ a point $y_n\in Y_n$ such that $g_{k(y_n)}\geq 1$ and $c_n(y_n)\to x$.
\end{prop}

\begin{thm}\label{rnsmumford}
Let $X$ be a Mumford curve over $\overline{\mbf Q}_p$. Then, $X$ satisfies resolution of non-singularities.
\end{thm}

\dem
Let $x\in X(\overline{\mbf Q}_p)$. According to proposition \ref{rnsprop}, it is enough to show that $x\in\overline V(X)$.

Let $\Omega=\mbf P^{1,\an}\backslash\mcal L$ be the topological universal cover of $X$, and let $\Gamma=\Gal(\Omega/X)\subset \PGL_2(\mbf C_p)$. Let $z$ be a point of $\Omega$ above $x$.
One can assume that $\Gamma\subset\PGL_2(K_0)$, where $K_0$ is a finite extension of $\mbf Q_p$. Let $g\in\Gamma$. Let $t\in \mcal L\cap \mbf P^1(\overline{\mbf Q}_p)$ be a point that is not fixed by $g$. Up to replacing $K_0$ by a finite extension, one can assume $t\in K_0$. Let $I=\{g^n(t)\}_{n\in\mbf Z}\subset K_0$.
According to proposition \ref{torsensinf}, there exists a finite subset $I_0$ of $I$ and a $\mbf Z_p(1)$-torsor $(c_n:Y_n\to \mbf P^1\backslash I_0)$ of $\mbf P^1\backslash I_0$ and a point $y_n$ of $Y_n$ such that $g_{k_{y_n}}\geq 1$ and $z_n:=c_n(y_n)\to z$. Let $x_n$ be the image of $z_n$ in $X$. Fix $n$ and show that $x_n\in\widetilde V(X)$. Let $\epsilon$ be small enough so that the canonical $\mu_{p^n}$ torsor of $\mbf G_m$ is split at $b_{1,\varepsilon}$. Let $f\in O(\mbf P^1\backslash I_0)^*$ such that $c_n=f^*\mu_{p^n}$.

Let $z_0\in\Omega(\mbf C_p)$. By replacing $f$ by $f(z)/f(z_0)$, one can
assume $f(z_0)=1$. Let $\Gamma'$ be a subgroup of $\Gamma$ of finite
index. Consider
$f_{\Gamma'}(z)=\prod_{g\in\Gamma'}\frac{f(g(z))}{f(g(z_0))}$: this product
converges uniformly on every affinoid domain of $\Omega$ and therefore
defines an element of $O(\Omega)$. Moreover, for every $g\in \Gamma'$,
$\frac{f_{\Gamma'}\circ g}{f_{\Gamma'}}$ is a constant function, \ie
$f_{\Gamma'}$ is a theta function of $X/\Gamma'$. If $(\Gamma_m)_{m\in\mbf
  N}$ is a decreasing sequence of subgroups of finite index of $\Gamma$
such that $\bigcap_m \Gamma_m=\{1\}$, the sequence
$(f_{\Gamma_m})_{m\in\mbf N}$ converges uniformly to $f$ on every affinoid
domain of $\Omega$. In particular, there exists a subgroup $\Gamma'$ of
finite index of $\Gamma$ such that
$|f_{\Gamma'}/f-1|_{z_n}<\varepsilon$. Let $c'=f_{\Gamma'}^*\mu_{p^n}:Y'\to
\Omega$. Then $c_n-c'$ is split at $z_n$. Therefore $c_n$ and $c'$ are
isomorphic above $z_n$ and there is $y'\in Y'$ above $z_n$ such that
$g_{k_{y'_n}}\geq 1$. Since $f_{\Gamma'}$ is a theta function of
$\Omega/\Gamma'$, there exists a $\mu_{p^n}$ torsor $c'':Y''\to
\Omega/\Gamma'$ such that $c'=p^*c''$ where $p$ is the topological cover
$\Omega\to\Omega/\Gamma'$. Since $Y'\to Y''$ is a topological cover, the
image $y''$ of $y'$ in $Y''$ is in $V(Y'')$. Since $Y''\to X$ is a finite
cover and maps $y''$ to $x_n$, $x_n\in\widetilde V(X)$. Therefore
$x\in\overline V(X)$.
\findem
Since $\P^1\backslash\{0,1,\infty\}$ and punctured Tate curves have finite
étale covers that are nonempty Zariski open subsets of Mumford curves, they
also satisfy resolution of nonsingularities. For every curve $X$ over
$\C_p$, there exists a nonempty Zariski open subset $U\subset X$ and a
finite cover $U\to\P^1\backslash\{x_1,\dots,x_n\}$ with $n\geq 3$. Since
$\P^1\backslash\{x_1,\dots,x_n\}$ satisfies resolution of non
singularities, $U$ also satisfy resolution of nonsingularities:
\begin{cor}
 Every $\C_p$-curve has a Zariski dense open subset which satisfies resolution of nonsingularities.
\end{cor}
\section{Resolution of non singularities and anabelian tempered geometry}
\subsection{Tempered fundamental group}
Let $K$ be a complete nonarchimedean field.

A morphism $f:S'\to S$ of $K$-analytic spaces is said to be an \emph{étale cover} if $S$ is covered by open subsets $U$ such that $f^{-1}(U)=\coprod V_j$ and $V_j\to U$ is étale finite~(\cite{dJ1}).

For example, finite étale covers, also called \emph{algebraic covers}, and covers in the usual topological sense for the Berkovich topology, also called \emph{topological covers}, are étale covers.

Then, André defines tempered covers as follows:
\begin{dfn} \label{def:rvt:temp}(\cite[def. 2.1.1]{andre1})
An étale cover $S' \to S$ is \emph{tempered} if it is a quotient of the composition of a topological cover $T' \to T$ and of a finite étale cover
 $T \to S$.
\end{dfn}
This is equivalent to say that it becomes a topological cover after pullback by some finite étale cover.

We denote by $\Covtemp(X)$ (resp. $\Covalg(X)$, $\Covtop(X)$) the category of tempered covers (resp. algebraic covers, topological covers) of $X$ (with the obvious morphisms).

\smallskip

A geometric point of a $K$-manifold $X$ is a morphism of Berkovich spaces $\mcal M(\Omega)\to X$ where $\Omega$ is an algebraically closed complete isometric extension of $K$.

Let $\bar x$ be a geometric point of $X$. Then one has a functor \[F_{\bar x}:\Covtemp(X)\to\Set\] which maps a cover $S\to X$ to the set $S_{\bar x}$. If $\bar x$ and $\bar x'$ are two geometric points, then $F_{\bar x}$ and $F_{\bar x'}$ are (non canonically) isomorphic (\cite[prop. 2.9]{dJ1}).

A functor $F:\Covtemp(X)\to\Set$ is said to be a \emph{fiber functor} if it is isomorphic to $F_{\bar x}$ for some (and therefore every) geometric point $\bar x$ of $X$. 

\begin{prop}
 A fiber functor of $\Covtemp(X)$ is pro-representable.
\end{prop}
If $F$ is a fiber functor of $\Covtemp(X)$, a pointed tempered cover of $X$ is a couple $(S,s)$ where $S$ is a tempered cover of $X$ and $s\in F(S)$.
A pro-tempered cover $\widetilde X$ of $X$ is called universal if $F_{\widetilde X}:=\Hom(\widetilde X,\ )$ is a fiber functor of $\Covtemp(X)$.

The tempered fundamental group of $X$ pointed at a fiber functor $F$ is \[\gtemp(X,F)=\Aut F.\]
The tempered fundametal group of $X$ pointed at a universal pro-tempered cover $\widetilde X$ is \[\gtemp(X,\widetilde X)=\Aut \widetilde X=\Aut F_{\widetilde X}.\]
The tempered fundamental group of $X$ pointed at a geometric point $\bar x$ is
\[\gtemp(X,\bar x)=\gtemp(X,F_{\bar x})\]
When $X$ is a smooth algebraic $K$-variety, $\Covtemp(X^{\an})$ and $\gtemp(X^{\an},F)$ will also be denoted simply by $\Covtemp(X)$ and $\gtemp(X,F)$.\\
By considering the stabilizers $(\Stab_{F(S)}(s))_{(S,s)}$ as a basis of open subgroups of $\gtemp(X,F)$, $\gtemp(X,F)$ becomes a topological group. It is a prodiscrete topological group.\\
When $X$ is algebraic, $K$ of characteristic zero and has only countably many finite extensions in a fixed algebraic closure $\overline K$, $\gtemp(X,F)$ has a countable fundamental system of neighborhood of $1$ and all its discrete quotient groups are finitely generated~(\cite[prop. 2.1.7]{andre1}).\\

Thus, as usual, the tempered fundamental group
depends on the basepoint only up to inner automorphism (this topological
group, considered up to conjugation, will sometimes be denoted simply by $\gtemp(X)$).\\
The full subcategory of tempered covers $S$ for which $F_{\bar x}(S)$ is finite is equivalent to $\Covalg(S)$, hence \[\widehat{\gtemp(X,\bar x)}=\ga(X,\bar x)\] (where $\widehat{\ }$ denotes here, and in the sequel, the profinite completion).\\
For any morphism $X\to Y$, the pullback defines a functor $\Covtemp(Y)\to\Covtemp(X)$. If $\bar x$ is a geometric point of $X$ with image $\bar y$ in $Y$, this gives rise to a continuous homomorphism \[\gtemp(X,\bar x)\to\gtemp(Y,\bar y)\] (hence an outer morphism $\gtemp(X)\to\gtemp(Y)$).\\
One has the analog of the usual Galois correspondence:
\begin{thm}\emph{(\cite[th. 1.4.5]{andre1})}\label{galcorr} A fiber functor $F$ induces an equivalence of categories between the category tempered covers of $X$ and the category $\gtemp(X,F)\text{-}\Set$ of discrete sets endowed with an action of $\gtemp(X,F)$ through factorizes through a finite quotient.\end{thm}

If $S$ is a finite Galois cover of $X$, its universal topological cover $S^{\infty}$ is still Galois and every connected tempered cover is dominated by such a Galois tempered cover.\\
Let $\bar x$ be a geometric point of $X$. Let $\text{``}\varprojlim_{I}\text{''} (S_i,s_i)$ be a pointed universal pro-finite cover of $X$. Let $(S_i^{\infty},s_i^\infty)$ be the pointed universal topological cover of $(S_i,s_i)$. Then $\text{``}\varprojlim_{I}\text{''} S_i^{\infty}$ is a universal pro-tempered cover of $X$.

If $S=\text{``}\varprojlim\text{''}S_i$ is a pro-tempered cover of $X$, one denotes by $|S|$ the topological space $\varprojlim S_i$.

Let $(\overline X,D)$ be a marked curve. Let $X=\overline X\backslash D$. If $S$ is a tempered cover of $X$, it extends uniquely, and functorially in $Y$, in a ramified cover $\overline S\to \overline X$ (\cite[th. III.2.1.11]{andre1}). If $S=\text{``}\varprojlim S_i\text{''}$ is a pro-tempered cover of $X$, one denotes by $|S|_c$ the topological space $\varprojlim \overline S_i$. Any morphism $S\to S'$ of pro-tempered covers induces 
a continuous map $|S|_c\to |S'|_c$. In particular, if $\widetilde X$ is a universal pro-tempered cover of $X$, $\gtemp(X,\widetilde X)=\Aut \widetilde X$ acts on $|\widetilde X|_c$.

\bigskip

Let $\mcal P$ be the set of prime numbers and Let $\mbb L$ be a subset of $\mcal P$. We call $\mbb L$-integer an integer which is a product of elements of $\mbb L$. One writes $(p')$ for $\mcal P\backslash\{p\}$.

A $\mbb L$-tempered cover $S$ of $X$ is a tempered cover such that there exists a finite étale Galois cover $Y\to X$ of index a $\mbb L$-integer such that $S\times_XY\to Y$ is a topological cover. We denote by $\Covtemp(X)^{\mbb L}$ the category of $\mbb L$-tempered covers. If $F$ is a fiber functor of $\Covtemp(X)$, we denote by $\gtemp(X,F)^{\mbb L}$ the topological group of automorphisms of $F_{|\Covtemp(X)^{\mbb L}}$. 
If $\widetilde X$ is a universal pro-tempered cover of $X$, one defines $\widetilde X^{\mbb L}=\text{``}\varprojlim_{\widetilde X\to Y}\text{''}Y$ where $Y$ runs over $\widetilde X$-pointed $\mbb L$-tempered cover.

If $\mbb L\subset \mbb L'$, the fully faithful functor $\Covtemp(X)^{\mbb L}\to\Covtemp(X)^{\mbb L'}$ induces a morphism $\gtemp(X,F)^{\mbb L'}\to\gtemp(X,F)^{\mbb L}$, and in particular when $\mbb L'=\mcal P$, one gets a morphism $\gtemp(X,F)\to\gtemp(X,F)^{\mbb L}$.

\subsection{Decomposition group of a point}
Let $(\overline X,D)$ be a marked $K$-curve, and let $X=\overline X\backslash D$.
Let $\widetilde X$ be a universal pro-tempered cover. One defines
$\widetilde X^{\mbb L}=\varprojlim_Y Y^\infty$  where $Y$ runs over the $\widetilde
X$-pointed finite étale covers of $X$ of index an $\mbb L$-integer and
$Y^\infty$ is the $\widetilde X$-pointed universal topological cover of
$Y$. If $p\notin \mbb L$, every morphism $Z\to Y$ of $\widetilde
X$-pointed finite étale covers of $X$ of index an $\mbb L$-integer
induces a morphism of tree $\mbb T_Z\to\mbb T_Y$; one defines $\mbb T^{\mbb L}$ to be $\varinjlim_{Y}\mbb T_Y$.
The family of embeddings $V(\mbb T_Y)\to Y^{\infty}$ induces an
embedding $V(\mbb T_X^{\mbb L})\to |\widetilde X^{\mbb L}|$. If $z=(z_Y)$ is an
edge of $\mbb T_X^{\mbb L}$, then every morphism $Y\to Y'$ induces a
homeomorphism $S_{z_Y}\to S_{z_Y'}$ and one defines $S_z=\varprojlim_Y
S_{z_Y}$.

Then $\gtemp(X,\widetilde X)$ acts on $|\widetilde X|_c$. Similarly $\gtemp(X,\widetilde X)^{\mbb L}$ acts on $|\widetilde X^{\mbb L}|_c$.
Let $x\in |\widetilde X|_c$. One denotes by $D_x$ the stabilizer of $x$ in $\gtemp(X,\widetilde X)$, and by $D_{x,\mbb L}$ its image in $\gtemp(X,\widetilde X)^{\mbb L}$. The group $D_{x,\mbb L}$ is the stabilizer of the image of $x$ in $|\widetilde X^{\mbb L}|_c$. The decomposition group depends only of the image of $x$ in $\overline X$ up to conjugacy.

\begin{thm}[{\cite[cor. 3.11]{mochi}}]\label{mochi0} If $X_{\alpha}$
  and $X_\beta$ are two hyperbolic $\overline{\Q}_p$-curves, every
  (outer) isomorphism $\gamma:\gtemp(X_{\alpha,\bb
    C_p})^{(p')}\simeq\gtemp(X_{\beta,\bb C_p})^{(p')}$ determines, functorially in $\gamma$, an isomorphism of graphs 
  $\bar\gamma:\mathbb G_{X_\alpha}\simeq\mathbb G_{X_\beta}$.\end{thm}
More precisely, the map $x\mapsto D_{x,(p')}$ identifies the vertices of $\mbb T_X$ with the maximal compact subgroups of $\gtemp(X)^{(p')}$ and two vertices $x,x'$ of $\mbb T_X$ are linked by an edge if and only if $D_{x,(p')}\cap D_{x',{(p')}}\neq \{1\}$. Therefore the isomorphism $\gamma:\gtemp(X_{\alpha,\bb
    C_p})^{(p')}\simeq\gtemp(X_{\beta,\bb C_p})^{(p')}$ induces an equivariant isomorphism of graphs $\mbb T_{X_\alpha}\to\mbb T_{X_\beta}$, which gives $\bar\gamma$ by quotienting by the action of the tempered fundamental group.

\begin{prop}\label{noncommens}Let $x_1\neq x_2\in \widetilde V(\widetilde X)$. Then $D_{x_1}$ and $D_{x_2}$ are not commensurable.\end{prop}
\dem Let $(Y,f:\widetilde X\to Y)$ be a pointed finite étale cover of $X$ such that $f(x_1)\neq f(x_2)\in V(Y)$. Then $D_{x_i}\cap \gtemp(Y,\widetilde X)=D_{f(x_1)}$. But the images of $D_{f(x_1)}$ and $D_{f(x_2)}$ in $\gtemp(Y,\widetilde X)^{(p')}$ are already not commensurable.
\findem

\begin{cor}\label{normalizer}
 Let $x\in\widetilde V(\widetilde X)$. Then $D_x$ is its own normalizer.
\end{cor}
\dem
Let $g$ be in the normalizer of $D_x$. Then $D_{g(x)}=gD_xg^{-1}=D_x$. Since $x,g(x)\in\widetilde V(\widetilde X)$, proposition \ref{noncommens} tells us that $x=g(x)$, \ie $g\in D_x$.
\findem

According to \cite[prop. 10]{gtem}, if $D$ is a compact subgroup of $\gtemp(X,\overline X)$ which is not a pro-$p$ group, there exists $x\in |\widetilde X|_c$ such that $D\subset D_x$.
Moreover, a point $x\in |\widetilde X|_c$ is in $\widetilde V(\widetilde X)$ if and only if there exists an open finite index subgroup $H\subset \gtemp(X)$ such that the image of $D_x\cap H$ in $H^{(p')}$ is non commutative.
Therefore, the set $\widetilde V(\widetilde X)$ can be identified with the set of conjugacy classes of maximal compact subgroups $D$ of $\gtemp(X)$ such that the image of $D\cap H$ in $H^{(p')}$ is non commutative for some open finite index subgroup $H$ of $\gtemp(X)$.

\subsection{Tempered theoreticness of Berkovich topology}
Let $(\overline X,D)$ be a $\overline Q_p$-marked curve and let $X=\overline X\backslash D$.
If $Y\to X$ is a Galois finite \'etale cover and $(\overline{\mcal Y},\mcal D_Y)$ is the stable model of $Y$, $\mcal X_Y:=\overline{\mcal Y}/\Gal(Y/X)$ is a semistable model of $\overline X$.

Let us say that a topological group is temp-like if it is isomorphic to the tempered fundamental group of a hyperbolic curve over $\overline Q_p$ that satisfies resolution of non-singularities. We will construct for any temp-like topological group $\Pi$ a topological space $\widetilde S(\Pi)$ endowed with a continuous action of $\Pi$. The construction will be purely group theoretic, so that it will be functorial with respect to isomorphism of topological groups. Moreover when $\Pi=\gtemp(X,\widetilde X)$, we will get a $\gtemp(X,\widetilde X)$-equivariant homeomorphism $|\widetilde X|\to\widetilde S(\gtemp(X,\widetilde X))$.

\bigskip

Let $\Pi$ be a temp-like topological group. We fix an isomorphism $\Pi\simeq\gtemp(X,\widetilde X)$ (but we will take care that the construction of $\widetilde S(\Pi)$ do not depend on this isomorphism). Then there is a smallest normal open subgroup $\Pi^{\infty}$ such that $\Pi/\Pi^\infty$ is torsionfree. One can also see $\Pi^\infty$ as the closed subgroup generated by the compact subgroups of $\Pi$. 
 One defines the topological group $\Pi^{(p')}:=\varprojlim_N \Pi/N^{\infty}$ where $N$ goes through open normal subgroups of finite index prime to $p$ of $\Pi$ (such an $N$ is also temp-like, so that $N^\infty$ is well defined).
The morphism $\Pi\to \Pi^{(p')}$ has dense image. We denote by $\Pi^{(p'),\infty}$ the kernel of $\Pi^{(p')}\to \Pi/\Pi^\infty$.

If $H$ is a normal open subgroup of $\Pi$ of finite index, let $\widetilde V(\Pi)_H$ be the set of maximal compact subgroups $D$ of $\Pi$ such that $D\cap H$ is not commutative. Let $\widetilde V(\Pi)=\bigcup_H \widetilde V(\Pi)_H$. The group $\Pi$ acts by conjugacy on $\widetilde V(\Pi)$ and on $\widetilde V(\Pi)_H$ for every $H$. 

Recall that there is an equivariant bijection $\widetilde V(\widetilde X)\to \widetilde V(\Pi)$ that maps $x$ to $D_x$. More precisely, If $(Y,f:\widetilde X\to Y)$ is the pointed finite étale cover of $X$ corresponding to $H$, it induces a bijection $\{x\in|\widetilde X|:f(x)\in V(Y)\}\to \widetilde V(\Pi)_H$. By quotienting by $H^\infty$, one gets a bijection $V(Y)^\infty\to \widetilde V(\Pi)_H/H^{\infty}$.

Let $H\subset \Pi$ be a normal open subgroup of finite index, $H$ is also temp-like: it is isomorphic to $\gtemp(Y,\widetilde X)$ for some Galois finite étale cover $Y$ of $X$.
Let $V(H)^{(p')}$ be the set of maximal compact subgroups of $H^{(p')}$. Let $E(H)^{(p')}$ be the set of pairs of elements $(D,D')$ of $V(H)^{(p')}$ such that $D\cap D'\neq \{1\}$. These data define a graph $\mbb G(H)^{(p')}$. Since $H$ is normal in $\Pi$, the group $\Pi$ acts on $H^{(p')}$ by conjugacy and therefore on $\mbb G(H)^{(p')}$.
Remark that $H^{(p')}$ also acts by conjugacy on $\mbb G(H)^{(p')}$ and that the action of $\Pi$ and of $H^{(p')}$ coincide on $H$.

According to \cite[Th. 6]{gtem}, there is a $\Pi$-equivariant isomorphism $\mbb T_Y^{(p')}\simeq \mbb G(H)^{(p')}$ that maps a vertex $x$ to its stabilizer $D_x$ by the action of $H^{(p')}$. If $e=(D,D')\in E(H)^{(p')}$, one denotes by $D_e:=D\cap D'\subset H^{(p')}$.
The group $D_e$ is the stabilizer of the image of $e$ in $\mbb T_Y^{(p')}$ for the action of $H^{(p')}$ on $\mbb G(H)^{(p')}$.

Let
\[\mbb G(H)=\mbb G(H)^{(p')}/H^{(p')}\quad \text{and} \quad \mbb G(H)^{\infty}=\mbb G(H)^{(p')}/H^{(p'),\infty}.\]
Then $\mbb G(H)$ can be identified with $\mbb G_Y$ and $\mbb G(H)^\infty$ with $\mbb T_Y$.

If $D\in \widetilde V(\Pi)_H$, the image of $D\cap H$ in $H^{(p')}$ is a maximal compact subgroup, and therefore defines an element of $V(H)^{(p')}$, hence a $\Pi$-equivariant map $p_H:\widetilde V(\Pi)_H\to V(H)^{(p')}$.

Moreover, the induced map $p^{\infty}_H:\widetilde V(\Pi)_H/H^{\infty}\to V(H)^{(p')}/H^{(p'),\infty}$ is bijective:
Indeed, the diagram
\[\xymatrix{\widetilde V(\Pi)_H/H^{\infty} \ar[r]\ar[rd] & V(H)^{(p')}/H^{(p'),\infty}\ar[d]\\
& V(Y^{\infty})
}
\]
is commutative and the two vertical maps are bijections.

One has $H^{\infty}\subset \Pi^{\infty}$ and one gets a $\Pi$-equivariant map $\iota_{H,\Pi}^{(p')}:H^{(p')}\to \Pi^{(p')}$.

\begin{lem}
Let $D\in V(H)^{(p')}$. The subgroup $\iota_{H,\Pi}^{(p')}(D)\subset \Pi^{(p')}$ is either:
\begin{itemize}
 \item an open subgroup  of $D_1$ for a unique $D_1\in V(\Pi)^{(p')}$;
\item an open subgroup of $D_e$ for a unique $e\in \overline E(\Pi)^{(p')}$;
\item $\{1\}$.
\end{itemize}
Let $e\in E(H)^{(p')}$. The subgroup $\iota_{H,\Pi}^{(p')}(D_e)\subset \Pi^{(p')}$ is either:
\begin{itemize}
\item an open subgroup of $D_e$ for a unique $e\in \overline E(\Pi)^{(p')}$;
\item $\{1\}$.
\end{itemize}
\end{lem}
\dem
Uniqueness is clear.

Since $p^{\infty}_H$ is surjective, up to conjugating $D$ by an element of $H^{(p')}$, one can assume that $D$ is the image of $p_H$. Let $D_0\in \widetilde V(\Pi)_H$ be a preimage of $D$ and let $x$ be the corresponding point of $\widetilde V(\widetilde X)$. Since $H\cap D_x$ is open in $D_x$, $\iota_{H,\Pi}^{(p')}(D)$ is open in the image $D_{x,(p')}$ of $D_x$ by $\Pi\to\Pi^{(p')}$. Let $\tilde x$ be the image of $x$ in $\widetilde X^{(p')}$. If $\tilde x\in V(\widetilde X^{(p')}$), then $D_{x,(p')}\in V(\Pi)^{(p')}$; if $\tilde x$ lies in an edge $e$ of $\mbb T^{(p')}_X$, then $D_{x,(p')}=D_e$; otherwise, the image of $x$ in $X$ lies in a disk and since every prime-to-$p$ cover of a disk is trivial, $D_{x,(p')}=\{1\}$.

if $e\in E(H)^{(p')}$, there exists $\bar y\in |\widetilde Y^{(p')}|$ such that $D_e=D_x$. Up to conjugating $D$ by an element of $H^{(p')}$, one can assume that there is $x\in |\widetilde X|$ which maps to $\bar y$. Let $\bar x$ be the image of $\bar y$ in $|\widetilde X^{(p')}|$. Once again, $\iota_{H,\Pi}^{(p')}(D_e)$ is open in $D_{\bar x}$. Since $\bar y\notin V(Y)^{(p')}$, $\bar x\notin V(X)^{(p')}$. Therefore either $\bar x$ lies in an edge $e'$ of $\mbb T^{(p')}_X$ and $D_{\bar x}=D_{e'}$ or $x$ lies outside the image of $S(X)^{(p')}$ and $D_{\bar x}=\{1\}$.
\findem

One denotes by $\mbb G_{\Pi^\infty}(H)=\mbb G(H)^{\infty}/\Pi^{\infty}$. If one identifies $\mbb G(H)^{\infty}$ with $\mbb T_Y$, then $\mbb G_{\Pi^\infty}(H)=\mbb T(\mcal X_Y)$.
The bijection $p_H^{\infty,-1}$ induces a bijection $V_{\Pi^{\infty}}(H)\to \widetilde V(\Pi)_H/\Pi^{\infty}$.
If $H'\subset H$ are two normal subgroups of $\Pi$, the $\Pi$-equivariant injective map $\widetilde V(\Pi)_H\to\widetilde V(\Pi)_H'$ induces an injective map $V_{\Pi^{\infty}}(H)\to V_{\Pi^{\infty}}(H')$.

Let $e=(D_1,D_2)\in E(\Pi)^{\infty}$. One denotes by 
\[\widetilde A^{(p')}_{e,H}:=\{D\in V(H)^{(p')}| \exists \tilde e\in E(\Pi)^{(p')},[\tilde e]=e\text{ and }\{1\}\neq\iota_{H,\Pi}^{(p')}(D)\subset D_{\tilde e}\}\]
The action of $H^{(p'),\infty}$ and of $\Pi^\infty$ on $V(H)^{(p')}$ stabilize $\widetilde A_{e,H}$.
 Let
\[A_{e,H}=(\widetilde A^{(p')}_{e,H}/H^{(p'),\infty})/\Pi^{\infty}\]
Thus $A_{e,H}$ is a subset of $V_{\Pi^\infty}(H)$. An element $D\subset \Pi$ of $\widetilde V(\Pi)_H$ is mapped to $A_{e,H}$ by $\widehat V(\Pi)_H\to V_{\Pi^\infty}$ if and only if the image of $D$ in $\Pi^{(p')}$ is a representative of $e$.
If $Y\to X$ is the pointed Galois cover corresponding to $H$,
then $A_{e,H}$ can be identified with $A_{z,Y}$, as defined in \eqref{eq:Amodel} where $z$ is the node of $\mcal X_Y$ corresponding to $e$. 

The full subgraph $\mbb G(A_{e,H})$ of $\mbb G_{\Pi^{\infty}}(H)$ with vertices $A_{e,H}\cup \{i_{H,\Pi}(D_1),i_{H,\Pi}(D_2)\}$ is a line (indeed the embedding $|\mbb G(A_{e,H})|\subset |\mbb T_{\mcal X_Y}|\subset X^{\infty}$ identifies $|\mbb G(A_{e,H})|$ with $\overline S_z\simeq [0,1]$ where $z$ is the edge of $\mbb T_X$ corresponding to $e$), so that $A_{e,H}$ is naturally a totally ordered set for which $D_1$ is the minimal element and $D_2$ is the maximal. If $\bar e=[(D_2,D_1)]\in E(\Pi)^\infty$ is the same edge of $\mbb G(\Pi)^\infty$ with the opposite orientation, then there is an obvious bijection $A_{e,H}\simeq A_{\bar e,H}$, which is decreasing.

Let $H'\subset H\subset \Pi$ be two finite index normal subgroup of $\Pi$, the injective map $V_{\Pi^{\infty}}(H)\to V_{\Pi^{\infty}}(H')$ maps $A_{e,H}$ to $A_{e,H'}$. The induced map $A_{e,H}\to A_{e,H'}$ is increasing. Let
\[A_e:=\varinjlim_H A_{e,H},\]
where $H$ goes through finite index normal subgroups of $\Pi$. By identifying $A_{e,H}$ with $A_{z,Y}$, one gets an increasing bijection $A_{e}\simeq A_{z}$. Since $X$ satisfies resolution of non-singularities, $A_{z}$ can be identified with the set of points of type $2$ of $S_z$. Thus $A_e$ is an ordered set which is, non-canonically, isomorphic to $\mbf Q\cap (0,1)$.
Let $\widehat A_e$ be the Dedekind completion of $A_e$: $\widehat A_e$ is an ordered topological space non-canonically isomorphic to $[0,1]$. The decreasing bijections $A_{e,H}\to A_{\bar e,H}$ are compatible and therefore induce a homeomorphism $\phi_e:\widehat A_e\to \widehat A_{\bar e}$
Let us call $0_e$ (resp. $1_e$) the minimal element of $\widehat A_e$. One then defines the topological space
\[S(\Pi)^{\infty}:=\left(V(\Pi)^{\infty}\coprod_{e\in E(\Pi)^{\infty}}\widehat A_e\right)/\sim\]
where $\sim$ is generated by \[\forall (D,D')\in E(\Pi)^\infty,D\sim 0_{(D,D')},\quad \forall e\in E(\Pi)^\infty\forall x\in\widehat A_e, x\sim \phi_e(x).\]
Since $\widehat A_e$ is non-canonically homeomorphic to $[0,1]$, $S(\Pi)^{\infty}$ is non-canonically homeomorphic to the geometric realization of $\mbb G(\Pi)^\infty$.

If $H$ is a finite index open normal subgroup of $\Pi$, one similarly gets a topological space $S(H)^\infty$ and the action of $\Pi$ on $H$ by conjugacy induces an action of $\Pi$ on $S(H)^\infty$.

Let $\tilde e\in E(H)^{(p')}$ and $\tilde e_0\in E(\Pi)^{(p')}$ be such that $\iota_{H,\Pi}^{(p')}(D_{\tilde e})$ is an open subgroup of $D_{\tilde e_0}$. Let $H'\subset H$ be an open normal subgroup of $\Pi$. One has $\widetilde A_{e,H'}\subset\widetilde A_{e_0,H'}$, as subsets of $V(H')^{(p')}$. Hence a map \[A_{e,H'}:=\widetilde A_{e,H'}/H\to A_{e_0,H'}:=A_{e_0,H'}.\] Since open normal subgroups of $\Pi$ which are inside $H$ are cofinal among open normal subgroups of $\Pi$ and among normal subgroups of $H$, one gets by taking colimits a map \[\alpha_{e,e_0}:A_e\to A_{e_0}\]

\begin{lem}\label{lemfunctsk}
There exists at most one continuous map $\psi_{H,\Pi}:S(H)^\infty \to S(\Pi)^\infty$ such that:
\begin{enumerate}[(i)]
 \item if $e\in E(H)^\infty$ and $\iota_{H,\Pi}^{(p')}(D_{\tilde e})=1$, then $\psi_{H,\Pi}$ is constant on $\widehat A_e$;
 \item if $\tilde e\in E(H)^{(p')}$ and $\tilde e_0\in E(\Pi)^{(p')}$ are such that $\iota_{H,\Pi}^{(p')}(D_{\tilde e})$ is an open subgroup of $D_{\tilde e_0}$, ${\psi_{H,\Pi}}_{|A_e}=\alpha_{e,e_0}$. 
\end{enumerate}
\end{lem}
\dem
Assume $\psi$ and $\psi'$ satisfy the condition. For every $\tilde e\in E(H)^{(p')}$ such that $\iota_{H,\Pi}^{(p')}(D_{\tilde e})\neq 1$, then thanks to $(ii)$, $\psi=\psi'$ on $A_e$ and therefore on $\widehat A_e$. For every $\tilde v\in V(H)^{(p')}$ such that $\iota_{H,\Pi}^{(p')}(D_{\tilde v})\neq 1$, there exists $\tilde e\in E(H)^{(p')}$ ending at $\tilde v$ such that $\iota_{H,\Pi}^{(p')}(D_{\tilde v})\neq 1$, and therefore $\psi(v)=\psi(v')$. Since $\mbb G(\Pi)^\infty$ is connected, one can link every edge and vertex of $\mbb G(\Pi)^\infty$ by a finite path to a vertex such that $\iota_{H,\Pi}^{(p')}(D_{\tilde v})\neq 1$. Up to reducing the path one can assume that for every edge $e$ of the path $\iota_{H,\Pi}^{(p')}(D_{\tilde e})=1$. By induction on the length of such a path, one gets that $\psi=\psi'$. using $(i)$.
\findem

\smallskip

There is a $\Pi$-equivariant injection $V(\Pi)^{(p')}\to X^{(p')}$ that maps $D$ to the unique $x\in X^{(p')}$ such that $D=D_x$. By quotienting by $\Pi^\infty$, one gets a $\Pi$-equivariant injection $V(\Pi)^{\infty}\to X^\infty$. Similarly, for every finite index subgroup $H$ of $\Pi$, there is a $\Pi$-equivariant injection $V(H)^\infty\to Y^\infty$, which induces by quotienting by $\Pi^\infty$ an injection  $V_{\Pi^\infty}(H)\to X^\infty$. If $H\subset H'$ are two finite index normal subgroups, the following diagram is commutative:
\[\xymatrix{V_{\Pi^\infty}(H)\ar[rd]\ar[r] & V_{\Pi^\infty}(H')
 \ar[d]\\
& X^{\infty}}
\]
Let $e$ be an edge of $\mbb G(\Pi)^\infty$, and let $z$ be the corresponding node of $\mcal X_k$. One gets compatible injective maps $A_{e,H}\to X^\infty$, whose image lie in $S_z$, hence an injective map $f_e:A_e\to S_z$, compatible with the reversing of edges. Since $A_e$ is dense in $\widehat A_e$, there is at most one extension of $f_e$ into a map $\widehat A_e\to X^{\infty}$, necessarily compatible with edges. 
Since $X$ satisfies resolution of non-singularities, the image of $A_e$ is exactly the set of points of type $(2)$ in $S_z$.
However, if $\mcal X$ is isomorphic to $\Spec O_K[X,Y]/(XY-a)$ in an étale neighborhood of $z$, then $\overline S_z$ can be identified with $[0,v(a)]$, and the set of points of type $(2)$ of $S_z$ is $\mbf Q\cap (0,v(a))$. On $A_{e,H}$ identified with a finite subset of $[0,v(a)]$, the tree structure is simply given by joining the consecutive points. Therefore $A_e\to [0,v(a)]$ is monotonous. Therefore $f_e$ extends to a unique homeomorphism $\tilde f_e:\widehat A_e\to S_z$, that preserves the end-points and is compatible with the reversing of edges. By gluing these maps, one therefore gets a continuous bijection $f^\infty:S(\Pi)^\infty\to S(X^\infty)$, which is a homeomorphism since $S(\Pi)^\infty$ is locally compact.

Similarly, if $H$ is a finite index open subgroup of $\Pi$, one gets a $\Pi$-equivariant homeomorphism $S(H)^\infty\to S(Y^\infty)$.

The composition $S(H)^\infty\simeq S(Y^{\infty})\stackrel{r_Xf^\infty\iota_Y}\to S(X^{\infty})\simeq S(\Pi)^\infty$, where $f^\infty$ is the map $Y^\infty\to X^\infty$, satisfies the properties of lemma \ref{lemfunctsk}. Therefore there exists a unique map $S(H)^\infty\simeq S(\Pi)^\infty$ satisying the properties of lemma \ref{lemfunctsk}.

If $H'\subset H$ are two finite index open subgroups of $\Pi$, then, since $H$ is also temp-like, there exists a unique map $S(H')^\infty\to S(H)^\infty$ satisfying the properties of lemma \ref{lemfunctsk} and this map is $\Pi$-equivariant by uniqueness. If $H''\subset H'$, then the diagram
\[
 \xymatrix{S(H'')^\infty \ar[rd] \ar[r] & S(H')^\infty\ar[d]\\ & S(H)^\infty}
\]is commutative. One therefore gets a projective system $(S(H)^{\infty})$ and one defines:
\[\widetilde S(\Pi)=\varprojlim_HS(H)^\infty.\]
The equivariant homeomorphisms induce a $\Pi$-equivariant homeomorphism $S(H)^\infty\to S(Y^\infty)$ induce a $\Pi$-equivariant homeomorphism \[\widetilde S(\Pi)\to\varprojlim_Y S(Y^\infty).\]
The maps $|\widetilde X|_c\to \overline Y^\infty\stackrel{r_Y}{\to} S(Y^\infty)$ are compatible and therefore induce a $\Pi$-equivariant map $|\widetilde X|_c\to \varprojlim_Y S(Y^\infty)$.
\begin{lem}
 The map $|\widetilde X|_c\to \varprojlim_Y S(Y^\infty)$ is a homeomorphism.
\end{lem}
\dem
First let us show that the map $r:\overline X^{\an}=|\widetilde X|_c/\Pi\to \varprojlim_YS(Y^\infty)/\Pi=\varprojlim_YS(\mcal X_Y)$ is a homemorphism (the proof is similar to the proof of prop. \ref{homeosst}).
Since $r_{\mcal X_Y}:\overline X^{\an}\to S(\mcal X_Y)$ is surjective for every $Y$ and $\overline X^{\an}$ is compact, $r$ is surjective.
Let $x_1\neq x_2\in \overline X^{\an}$ and show that $r(x_1)\neq r(x_2)$. One can assume $r_X(x_1)=r_X(x_2)$. Let $[x_1,x_2]$ be the smallest subset of $X$ containing $x_1$ and $x_2$, endowed with the total order such that $x_1<x_2$. Let $y_1<y_2$ be two points of type 2 in $[x_1,x_2]$ and let $Y\to X$ be a finite Galois cover such that $y_1,y_2\in V(\mcal X_Y)$. Then $r_{\mcal X_Y}(x_1)<y_1<y_2<r_{\mcal X_Y}(x_2)$ in $[x_1,x_2]$, which proves the injectivity of $r_{\mcal X_Y}$ and therefore of $r$. Since $\overline X^{\an}$ is compact, $r$ is an homeomorphism.

By pulling back along $S(X^\infty)\to S(X)$, one gets that \[r^\infty:\overline X^\infty=\overline X^{\an}\times_{S(X)}S(X^\infty)\to\varprojlim_Y S(\mcal X_Y)\times_{S(X)}S(X^\infty)=\varprojlim_Y S(\mcal X_Y^{\infty})\] is a homeomorphism.

Similarly one gets that the map $\overline Y^\infty\to \varprojlim_Z S(\mcal Y_Z^{\infty})$, where $Z$ goes through Galois pointed cover of $Y$ (one can even restrict to $Z$ over $Y$ Galois over $X$ since they are cofinal among Galois pointed cover of $Y$) is a homeomorphism for every $Y$ Galois.
Therefore the map \[|\widetilde X|_c\to\varprojlim_{Z\to Y\to X} S(\mcal Y_Z^\infty)\to\varprojlim_Z S(Z^\infty)\] is a homeomorphism (the right arrow is a homeomorphism because the full subcategory of the category  of morphisms $Y\to Z$ of pointed Galois cover over $X$ which consists of isomorphisms is a cofinal category).
\findem
One thus gets an equivariant homeomorphism
\begin{equation}\label{isomunivcov}
 \widetilde S(\Pi)\to |\widetilde X|_c.
\end{equation}

\bigskip

Let $(\overline X,D)$ be a $\overline Q_p$-marked curve and let $X=\overline X\backslash D$.
If $Y\to X$ is a Galois finite \'etale cover and $(\overline{\mcal Y},\mcal D_Y)$ is the stable model of $Y$, $\mcal X_Y:=\overline{\mcal Y}/\Gal(Y/X)$ is a semistable model of $\overline X$ and one gets a refinement $\phi_Y:\mbb G_{Y/X}:=\mbb G_Y/\Gal(Y/X)=\mbb G_{\mcal Y_s/\Gal(Y/X)}\to\mbb G_X$.

If $X$ satisfies non resolution of singularities, the family $(\overline{\mcal Y}/\Gal(Y/X))_Y$ is cofinal among semistable models of $\overline X$. Thus, if $e\in \mcal E(X)$, $A_e=\varinjlim A_{\phi_Y,e}$ and $\overline X^{\an}\to\varprojlim_{Y} |\mbb G_{Y/X}|_{\can}$ is a homeomorphism.

\begin{thm}\label{resnonsing}
 Let $X_1=\overline X_1\backslash D_1,X_2=\overline X_2\backslash D_2$ be two marked curves satisfying resolution of non singularities and let $\widetilde X_i$ be a universal pro-tempered cover of $X_i$. Let $\psi:\gtemp(X_1,\widetilde X_1)\simeq\gtemp(X_2,\widetilde X_2)$ be an isomorphism. Then there exists a unique homeomorphism $\bar\psi:|\widetilde X_1|_c\to|\widetilde X_2|_c$ which is $\gtemp$-equivariant in the sense that the following diagram commutes:
\[\xymatrix{\gtemp(X_1,\widetilde X_1)\times |\widetilde X_1|_c \ar[r]\ar[d]^{\psi\times\bar\psi} & |\widetilde{X}_1|_c\ar[d]^{\bar\psi}\\
\gtemp(X_2,\widetilde X_2)\times |\widetilde X_2|_c \ar[r] & |\widetilde X_2|_c}
\]
In particular, by quotienting by the tempered fundmental group, one gets a homeomorphism:
\[\overline X^{\an}_1\to\overline X_2^{\an}\]

\end{thm}

\dem
First, assume $\bar\psi_1,\bar\psi_2:|\widetilde X_1|_c\to|\widetilde X_2|_c$ are two $\gtemp$-equivariant homeomorphisms. Then, if $x\in |\widetilde X_1|_c$, $\psi(D_{x})=D_{\bar\psi_1(x)}=D_{\bar\psi_2(x)}$. If $x$ is in $\widetilde V(\widetilde X_1)$, then there exists an open subgroup $H$ of finite index of $\gtemp(X_1,\widetilde X_1)$ such that $(D_x\cap H)^{(p')}$ is not commutative. Then $\psi(H)$ is a finite index subgroup of $\gtemp(X_2,\widetilde X_2)$ and $(D_{\bar\psi_1(x)}\cap \psi(H))^{(p')}$ and $(D_{\bar\psi_2(x)}\cap \psi(H))^{(p')}$ are not commutative. Therefore $\bar\psi_1(x)$ and $\bar\psi_2(x)$ are of type 2 and have the same decomposition group: according to proposition \ref{noncommens}, $\bar\psi_1(x)=\bar\psi_2(x)$ for every point of $\widetilde V(\widetilde X_1)$. Since $\widetilde V(\widetilde X_1)$ is dense in $|\widetilde X_1|_c$, one has $\bar \psi_1=\bar\psi_2$.

The morphism $\psi$ induces an equivariant homeomorphism $\widetilde S(\gtemp(X_1,\widetilde X_1))\to \widetilde S(\gtemp(X_2,\widetilde X_2))$. One gets from \eqref{isomunivcov} equivariant homeomorphisms
$|\widetilde X_1|_c\to \widetilde S(\gtemp(X_1,\widetilde X_1))$ and $|\widetilde X_2|_c\to \widetilde S(\gtemp(X_2,\widetilde X_2))$. One gets the wanted isomorphism by composition.

\findem

\begin{prop}
 Let $\psi_a,\psi_b:\gtemp(X_1,\widetilde X_1)\to\gtemp(X_2,\widetilde X_2)$ be two isomorphisms. If $\tilde \psi_a=\tilde\psi_b$, then $\psi_a=\psi_b$.
\end{prop}
\dem
Let $g\in\gtemp(X_1,\widetilde X_1)$.
If $x_1\in|\widetilde X_1|$, then
\[\psi_a(g)D_{\tilde\psi_a(x_1)}\psi_a(g)^{-1}=D_{\tilde\psi_a(gx_1)}=D_{\tilde\psi_b(x_1)}= \psi_b(g)D_{\tilde\psi_b(x_1)}\psi_b(g)^{-1}=\psi_b(g)D_{\tilde\psi_a(x_1)}\psi_b(g)^{-1}\]
which implies that $g_0:=\psi_a(g)^{-1}\psi_b(g)$ is in the normalizer $N_{\psi_a(x_1)}$ of $D_{\psi_a(x_1)}$. Since $\psi_a$, is bijective, $g_0\in\bigcap_{x_2\in\widetilde X_2}N_{x_2}$. If $x_2\in \widetilde V(\widetilde X_2)$, then $N_{x_2}=D_{x_2}$. Therefore $g_0(x_2)=x_2$ for every $x_2\in\widetilde V(\widetilde X_2)$. Since $X_2$ satisfies resolution of non singularities, $\widetilde V(\widetilde X_2)$ is dense in $\widetilde X_2$, and thus $g_0(x)=x$ for every $x\in\widetilde X_2$. If $x$ is of type 1, then $D_x=\{1\}$. Therefore $g_0=1$, \ie $\psi_a(g)=\psi_b(g)$.
\findem

In particular, if $(X,\widetilde X)$ is a pointed curve satisfying resolution of non singularities, one has an injective morphism of groups: $\Aut\gtemp(X,\widetilde X)\to \Aut |\widetilde X|$.

\section{Tate curves}
Let $q_1,q_2\in\overline{\mbf Q_p}$ such that $|q_1|<1$ and $|q_2|<1$. Let $E_i=\Gm/q_i^{\mbf Z}$ and let $X_i=E_i\backslash\{1\}$. 

If there exists $\sigma\in G_{\mbf Q_p}$ such that $q_1=\sigma(q_2)$, there is a $\mbf C_p$-isomorphism $X_1\simeq X_2\otimes_{\mbf C_p}\mbf C_p$ where the base change $\mbf C_p\to\mbf C_p$ is $\sigma$. Therefore $X_1$ and $X_2$ are isomorphic analytic spaces over $\mbf Q_P$ and therefore have isomorphic tempered fundamental group.
The following theorem states that the converse is also true:

\begin{thm}\label{tatecurves}
Let $\psi:\gtemp(X_1,\widetilde X_1)\simeq\gtemp(X_2,\widetilde X_2)$ be an isomorphism.
 There exists $\sigma\in G_{\mbf Q_p}$ such that $q_1=\sigma(q_2)$, \ie $E_1$ and $E_2$ are isomorphic analytic spaces over $\Q_p$.
\end{thm}
\begin{rem}
The curves $X_1$ and $X_2$ satisfy the assumptions of theorem \ref{resnonsing}, and thus $\psi$ induces a homeomorphism $E_1^{\an}\to E_2^{\an}$. However, the author does not know, even in this situation, if this homeomorphism comes from an analytic morphism.
 The author does not know how to associate to $\psi$ a particular $\sigma$.
\end{rem}
\dem
Recall that $|q_1|=|q_2|$ and that $\psi$ induces a unique equivariant homeomorphism $\tilde\psi:|\widetilde X_1|_c\simeq |\widetilde X_2|_c$. The induced homeomorphism $E_1\simeq E_2$ maps $X_1$ onto $X_2$
Let $\Omega_i$ be the universal $\widetilde X_i$-pointed topological cover of $X_i$ and let $\bar\psi: |\overline{\Omega}_1|\to|\overline{\Omega}_2|$ the homeomorphism induced by $\tilde \psi$. 

Let $E_{i,l}$ be the unique $\widetilde X_i$-pointed connected topological cover of $E_i$ of degree $l$. Let $X_{i,l}=X_i\times_{E_i}E_{i,l}$. The isomorphism $\psi$ induces an isomorphism $\psi_l:\gtemp(X_{1,l},\widetilde X_1)\to \gtemp(X_{2,l},\widetilde X_2)$, whence an isomorphism \[\psi_{l,n}:H^1(X_{2,l},\mu_n)\to H^1(X_{1,l},\mu_n)\] functorial in $l$ for divisibility.

Let $\mbb G_i$ be the semigraph of $X_i$ and let $\mbb T_i$ be its universal cover. Let $g$ be the isomorphism $\mbb T_1\to\mbb T_2$ induced by $\psi$. One identifies $\overline{\Omega}_i$ with $\Gm$ in such a way that $\bar\psi(1)=1$ and $\bar\psi(q_1)=q_2$. If $j\in \Z$, one denotes by $e_{i,j}$ the cuspidal edge of $\mbb T_i$ corresponding to $q_i^j\in\overline{\Omega}_i\backslash\Omega_i$ and by $v_{i,j}$ the vertex of $\mbb T_i$ at which $e_{i,j}$ ends. Since $\bar\psi(1)=1$ and $\bar\psi(q_1)=q_2$, $g(e_{1,j})=e_{2,j}$ and $g(v_{1,j})=v_{2,j}$. One denotes by $e'_{i,j}$ the unique oriented edge joining $v_{i,j}$ to $v_{i,j+1}$.

All the cohomology groups will be cohomology groups for \'etale cohomology in the sense of algebraic geometry or in the sense of Berkovich. (one can replace  \'etale cohomology of $X^{\an}$ by \'etale cohomology of $X$ thanks to \cite[thm. 3.1]{berkgaga}).
Since $\Gamma_l\simeq \Z$, one has $H^2(\Gamma,\mu_n)=0$. Therefore the spectral sequence 
\[
H^p(\Gamma_l,H^q(\Omega,\mu_n))\implies H^{n}(X_{i,l},\mu_n)
\]
of the Galois étale cover $\Omega_i\to X_{i,l}$ gives us an exact sequence of cohomology groups for Berkovich étale topology:
\[1\to \Hom(\Gamma_l,\mu_n)\to H^1(X_{i,l},\mu_n)\to H^1(\Omega_i,\mu_n)^{\Gamma_l}\to 1.\]
The map $O^*(\Omega_i)/O^*(\Omega_i^n)\to H^1(\Omega_i,\mu_n)$ given by Kummer theory (see
\cite[prop. 4.1.7]{berketale} for the Kummer exact sequence in Berkovich
\'etale topology) induces an morphism
\[
\big(O^*(\Omega_i)/O^*(\Omega_i)^n\big)^{\Gamma_l}\to H^1(\Omega_i,\mu_n)^{\Gamma_l},
\]
which turns to be an isomorphism (cf. \cite[\S 1.4.2]{pia});
hence an exact sequence
 \[1\to \Hom(\Gamma_l,\mu_n)\to H^1(X_{i,l},\mu_n)\to \big(O^*(\Omega_i)/O^*(\Omega_i)^n\big)^{\Gamma_l}\to 1.\]
where $\Gamma_l=\Gal(\Omega_i/X_{i,l})$.  

One can describe $O^*(\Omega_i)$ in terms of currents (as done in \cite{vdP1} for Mumford curves).
If $A$ is a ring, a $A$-current on $\mbb T_i$ is a function $c:\{e_{i,j}\}_{j\in\Z}\coprod \{e'_{i,j}\}_{j\in\Z}\to A$ such that for every $j\in\Z$, $c(e'_{i,j+1})=c(e'_{i,j})+c(e_{i,j+1})$. Let $C(\T_i,A)$ be the $A$-module of $A$-currents on $\mbb T_i$.
There is a natural isomorphism $\alpha_i:\C(\T_i,\Z)\to O^*(\Omega_i)/\C_p^*$ defined by \[\alpha_i(c)=x^{c(e'_{i,0})}\prod_{j\geq 1}\big(\frac{x-q_i^j}{x}\big)^{c(e_{i,j})}\prod_{j\leq 0}\big(\frac{x-q_i^j}{q_i^j}\big)^{c(e_{i,j})}\]
Conversely, if $f\in O^*(\Omega_i)$, one can compute $\alpha_i^{-1}(f)$ in the following way. For every $j\in \Z$, the restriction of $f$ to the open annulus $U_j=\{z\in\P^{1,\an}, |q_i^j|<|x|_z<|q_i^{j+1}|\}$ can be written in a unique way as $f(x)=x^{m_j}g_j(x)$ where $m_j\in\Z$ and $|g_j|$ is constant on $U_j$. One has $\alpha_i^{-1}(f)(e'_{i,j})=m_j$. Similarly, the restriction of $f$ to the punctured open disk $V_j=\{z\in\P^{1,\an}, 0<|x-q_i^j|_z<|q_i^j|\}$ can be written in a unique way as $f(x)=x^{n_j}f_j(x)$ where $n_j\in\Z$ and $|f_j|$ is constant on $V_j$. One has $\alpha_i^{-1}(f)(e_{i,j})=n_j$.

Therefore, one gets an isomorphism $\alpha_{i,n}:O^*(\Omega_i)/O^*(\Omega_i)^n\to C(\T_i,\Z/n\Z)$, hence an exact sequence:
\[1\to \Hom(\Gamma_l,\mu_n)\to H^1(X_{i,l},\mu_n)\to C(\T_i,\Z/n\Z)^{\Gamma_l}\to 1.\]

Since $\varinjlim_l \Hom(\Gamma_l,\mu_n)=0$, it induces an isomorphism \[a_i:\varinjlim_l H^1(X_{i,l},\mu_n)\to C(\mbb T_i,\Z/n\Z)^{(\Gamma)},\] where $C(\mbb T_i,\mbb Z/n\mbb Z)^{(\Gamma)}$ is the set of $\mbb Z/n\mbb Z$-currents on $\mbb T_i$ that are invariant under some finite index subgroup of $\Gamma$.

Consider the diagram:
\[\xymatrix{\varinjlim_l H^1(X_{2,l},\mu_n) \ar[r]^{f:=\varinjlim_l\psi_{l,n}}\ar[d]^{a_1}  & \varinjlim_l H^1(X_{1,l},\mu_n)  \ar[d]^{a_2} \\ C(\mbb
  T_2,\mbb Z/n\mbb Z)^{(\Gamma)} \ar[r]^{g^*} & C(\mbb T_1,\mbb Z/n\mbb Z)^{(\Gamma)} }
\]
where the lower arrow is induced by $g:\mbb T_1\to\mbb T_2$.

The following lemma shows that this diagram is commutative up to a constant in $(\Z/n\Z)^*$ (cf. \cite[prop. 13]{pia} for a similar result for Mumford curves of genus greater than 2).
\begin{prop}\label{compatibilitycurrents}There exists a unique $\alpha\in(\mbf Z/n\mbf Z)^*$ such that $a_2f=\alpha g^*a_1$.
\end{prop}
\dem
Let $\tilde f=a_2fa_1^{-1}$. We have to show that there exists
$\lambda\in(\mbf Z/n\mbf Z)^*$ such that for every $c\in C(\mbb
  T_2,\mbb Z/n\mbb Z)^{(\Gamma)}$ and every edge $e$ of $\mbb T_1$, $\tilde f(c)(e)=\lambda c(g(e))$.\\
Let $j\in \Z$. According to \cite[lem. 4.2]{metricmumford}, a finite cover of $X_{1,l}$ is ramified at $e_{1,j}$ if and only if the corresponding cover of $X_{2,l}$ is unramified at $e_{2,j}$. Therefore $\tilde f(c)(e_{1,j})=0$ if and only if $c(e_{2,j})=0$. Therefore there exists $\lambda_j\in(\mbf Z/n\mbf Z)^*$ such that $\tilde f(c)(e_{1,j})=\lambda_j c(e_{2,j})$ for every $c\in C(\mbb
  T_2,\mbb Z/n\mbb Z)^{(\Gamma)}$. Let $l$ be a positive integer. Let $j\in \Z$ and let $c_j$ be the  $\Gamma_{l}$-invariant current defined by
\[c_j(e_{1,k})=\left\{\begin{array}{ll}
                      1 & \text{if}\quad k=j \mod l\\
-1 &\text{if}\quad k=j+1 \mod l\\
0 &\text{otherwise}
                     \end{array}\right.
\]
\[c_j(e'_{1,k})=\left\{\begin{array}{ll}
                      0 & \text{if}\quad k\in [j+1,j+l-1] \mod l\\
1 &\text{if}\quad k=j \mod l
                     \end{array}\right.
\]
If $l$ is big enough (for example if $l\geq 2+2\frac{v_p(n)+2}{|q_1|}$ according to \cite[cor. 4.10]{metricmumford}), the $\mu_n$-torsor corresponding to $c_j$ is split at  $v_{1,j+\lceil\frac{l+1}{2}\rceil}$. Therefore, the $\mu_n$-torsor corresponding to $\tilde f(c_e)$ is split at  $v_{2,j+\lceil\frac{l+1}{2}\rceil}$, which implies, according to \cite[prop. 4.11]{metricmumford}, that $\tilde f(c_j)$ is zero at all the edges ending at $v_{2,j+\lceil\frac{l+1}{2}\rceil}$. Then $\tilde f(c_j)(e_{2,j})=\lambda_j$, $\tilde f(c_j)(e_{2,j+1})=\lambda_{j+1}$, $\tilde f(c_j)(e_{2,k})=0$ for all $k\neq j,j+1 \mod l$,  and is zero for some non cuspidal edge between $v_{2,j+1}$ and $v_{2,j+l}$. Therefore $\tilde f(c_j)=\lambda_jg^*(c_j)$ and $\lambda_j=\lambda_{j+1}$. One thus gets that $\lambda_j$ does not depend on $j$, one simply denotes it by $\lambda$.  The group of $\Gamma_{l}$-equivariant current is generated by $(c_j)_{j\in [0,2l-1]}$ so that one gets $\tilde f(c)=\lambda g^{*}(c)$ for every current $c$, which ends the proof.

\findem

Let $A_i$ be the multiplicative group of non-zero meromorphic functions on $\Gm$ with no poles and no zeroes on $\Omega_i\subset \Gm$. Let $B_i\subset A_i$ be the subgroup of $A_i$ consisting of functions for which $1$ is neither a pole nor a zero. Let $A'_i\subset\Omega^1(\Omega_i)$ be the groupe of regular differentials on $\Gm$ with no poles on $\Omega_i$.

The map $d\log\circ\alpha_i: C(\T_i,\Z)\to A'_i$ can be extended by linearity to a map 
$\delta_i:C(\T_i,\Z_p)\to A'_i$ defined by
\[\delta_i(c)=c(e'_{i,0})\frac{dx}{x}+\sum_{j\geq 1}c(e_{i,j})(\frac{dx}{x-q_i^j}-\frac{dx}{x})+\sum_{j\leq 0}c(e_{i,j})\frac{dx}{x-q_i^j}.\]
If $z\in\Gm(\C_p)$, one denotes by and $\omega\in A'_i$, one denotes by $\ord_z(\omega)\in \Z_{\geq -1}$ the $(x-z)$-adic valuation of $\frac{\omega}{dx}\in \frac{1}{x-z}\C_p[[x-z]]$.

\begin{lem}\label{lemordre}
 Let $c$ be in $C(\T_1,\Z_p)$, and let $z\in\Gm(\C_p)$. Assume that $\bar\psi(z)$ is also of type $1$, \ie $\bar\psi(z)\in\Gm(\C_p)$. Then $\ord_z(\delta_1(c))=\ord_{\bar\psi(z)}(\delta_2(g^*(c)))$.
\end{lem}
\dem
Let $J=\{j\in\Z|C(e_{1,j})\neq 0\}$. For $z=q_1^j$ with $j\in J$, then $\bar\psi(q_1^j)=q_2^j$ and thus $\ord_{q_1^j}(\delta_1(c))=\ord_{q_2^j}(\delta_2(g^*(c)))=-1$. One can thus assume $z\in\Gm(\C_p)\backslash\{q_1^j\}_{j\in J}$ and therefore $\ord_z(\delta_1(c)),\ord_{\bar\psi(z)}(\delta_2(g^*(c)))\geq 0$.

Let $c_n\in C(\T_1,\Z)^{(\Gamma)}$ be such that $|c_n(e_{1,j})-c(e_{1,j})|,|c_n(e'_{1,j})-c(e'_{1,j})|\leq p^{-n}$ if $-|j|_{\infty}\log_p |q_1|\leq n+\frac{p}{p-1}$.

Then $\delta_1(c_n)\to \delta_1(c)$ on every affinoid subspace of $\Gm\backslash\{q_1^j\}_{j\in J}$. Similarly 
$\delta_2(g^*(c_n))\to \delta_2(g^*(c))$ on every affinoid subspace of $\Gm\backslash\{q_2^j\}_{j\in J}$.
Therefore
\[\ord_z(\delta_1(c))\geq m\iff \forall \epsilon>0\ \exists N\ \forall n\geq N, \sum_{z'\in D(z,\epsilon)}\ord_{z'}(\delta_1(c_n))\geq m.\]
Therefore it is enough to prove the lemma for every $c_n$. One thus assumes that $c\in C(\T_1,\Z)^{(\Gamma)}$. Let $N$ be such that $c\in C(\T_1,\Z)^{\Gamma_N}$. The isomorphism $\psi$ induces an isomorphism $\psi_N:\gtemp(X_{1,N})\to\gtemp(X_{2,N})$. 

Consider the image $Y_1$ of $c$ under the map
\[C(\T_1,\Z)\to O^*(\Omega_1)\to H^1(\Omega_1,\Z_p(1)),\]
and let $Y_{1,n}$ be the induced $\mu_{p^n}$-torsor on $\Omega_1$; it extends in an unramified $\mu_{p^n}$-torsor of $\Omega_1\backslash\{q_1^j\}_{j\in J}$. Let $z_{1,n}$ be the point of $\Omega_1$ such that $Y_{1,n}$ is not split at $z_{1,n}$ but is split on $[z,z_{1,n})$.
Consider the image $Y_2$ of $g^*(c)$ under the map $C(\T_2,\Z)\to O^*(\Omega_2)\to H^1(\Omega_2,\Z_p(1))$, let $Y_{2,n}$ be the induced $\mu_{p^n}$-torsor on $\Omega_2$ and let $z_{2,n}$ be the point of $\Omega_2$ such that $Y_{2,n}$ is not split at $z_{2,n}$ but is split on $[\bar\psi(z),z_{2,n})$. Since $c$ is $\Gamma_N$-invariant, $Y_{i,n}$ is in the image of $\Hom(\gtemp(X_{i,N}),\mu_{p^n})=H^1(X_{i,N},\mu_{p^n})\to H^1(\Omega_i,\mu_{p^n})$. For every preimage $\beta_i$ of $Y_{i,n}$ in $\Hom(\gtemp(X_{i,N}),\mu_{p^n})$, for every $Y_{i,n}$ is split at a point $z'\in\Gm$ if and only if $D_{z'}\subset\Ker(\beta_i)$, where $D_{z'}$ is a decomposition group of $z'$ in $\gtemp(X_{i,N})$. According to proposition \ref{compatibilitycurrents}, there exists $\alpha\in (\Z/p^n\Z)^*$ such that, if $\beta_1$ is a preimage of $Y_{1,n}$ in $\Hom(\gtemp(X_{1,N}),\mu_{p^n})$, then $\alpha\beta_1\psi^{-1}_N$ is a preimage of $Y_{2,n}$ in $\Hom(\gtemp(X_{2,N}),\mu_{p^n})$. Since $\Ker(\alpha\beta_1\psi^{-1}_N)=\psi(\Ker(\beta_1))$ and $D_{\bar\psi(z')}=\psi_N(D_{z'})$, one gets that $Y_{1,N}$ is split at $z'$ if and only if $Y_{2,N}$ is split at $\bar\psi(z')$, and therefore $z_{2,n}=\bar\psi(z_{1,n})$.

\smallskip

Let $c_0\in C(\T_1,\Z)^\Gamma$ be defined by $c_0(e_{1,i})=0$ and $c_0(e'_{1,i})=1$ for every $i\in\Z$. Then $\alpha_1(c_0)(x)=x$ and $\alpha_2(g^*c_0)(x)=x$. Let $\phi_{1,n}$ (resp. $\phi_{2,n}$) be a preimage of $c_0\mod p^n$ (resp. $g^*c_0 \mod p^n$) by the map $\Hom(\gtemp(X_1,\widetilde X_1),\mu_{p^n})=H^1(X_1,\mu_{p^n})\to C(\T_1,\mu_{p^n})^\Gamma$ (resp. $\Hom(\gtemp(X_2,\widetilde X_2),\mu_{p^n})=H^1(X_2,\mu_{p^n})\to C(\T_2,\mu_{p^n})^\Gamma$).
Let $z'_{1,n}=b_{z,|z|p^{-n-\frac{1}{p-1}}}\in \Gm$ and $z'_{2,n}=b_{\bar\psi(z),|\bar\psi(z)|p^{-n-\frac{1}{p-1}}}\in \Gm$. According to \cite[lem. 4.2]{metricmumford}, $z'_{1,n}$ is characterized in $\Gm$ by the fact that $c_{\can,p^n}$ is not split at $z'_{1,n}$ but is split above $[z,z'_{1,n})$.
Therefore $z'_{1,n}$ is also characterized by the fact that $D_{z'_{1,n}}\nsubseteq \Ker\phi_{1,n}$ and $D_{z'}\subset \Ker\phi_{1,n}$ for every $z'\in [z,z'_{1,n})$. Similarly, $z'_{2,n}$ is characterized by the fact that $D_{z'_{2,n}}\nsubseteq \Ker\phi_{2,n}$ and $D_{z'}\subset \Ker\phi_{2,n}$ for every $z'\in [z,z'_{2,n})$.
According to \ref{compatibilitycurrents}, one can choose $\phi_{2,n}$ to be $\alpha\phi_{1,n}\psi^{-1}$, so that $\Ker\phi_{2,n}=\psi(\Ker\phi_{1,n})$. Therefore, since $\bar\psi$ is compatible with decomposition groups, $z'_{2,n}=\bar\psi(z'_{1,n})$.

Using proposition \ref{rayonconv}, one gets
 \begin{align*}\ord_z(\delta_1(c))+1 & =\lim_n\frac{1}{n}\inf\{m,z'_{1,m}\in [z,z_{1,n}]\}\\
& = \lim_n\frac{1}{n}\inf\{m,\bar\psi(z'_{1,m})\in [\bar\psi(z),\bar\psi(z_{1,n})]\}\\
& = \lim_n\frac{1}{n}\inf\{m,z'_{2,m}\in [\bar\psi(z),z_{2,n}]\}\\
& = \ord_{\bar\psi(z)}(\delta_2(g^*c))+1.\end{align*}
\findem

Let $\mu:\N_{>0}\to \{-1,0,1\}$ be the Moebius function. Let $n\geq 1$. Let $c_n\in C(\T_1,\Z)$ defined by
\begin{itemize}
 \item $c_n(e_{1,j})=0$ if $j\leq 0$;
 \item $c_n(e'_{1,j})=0$ if $j\leq 0$;
 \item $c_n(e_{1,j})=\mu(\frac{j}{n})$ if $j\geq 1$ and $j=0\mod n$;
 \item $c_n(e_{1,j})=0$ if $j\geq 1$ and $j\neq 0\mod n$;
 \item $c_n(e'_{1,j})=\sum_{k=1}^{\lfloor \frac{j}{n}\rfloor}\mu(k)$ if $j\geq 1$.
\end{itemize}

The associated differentials are: 
\[\begin{array}{rcl}
\delta_1(c_n) & = & \sum_{j\geq 1}\mu(j)\big(\frac{1}{x-q_1^{nj}}-1\big)dx\\
\delta_2(g^*(c_n)) & = & \sum_{j\geq 1}\mu(j)\big(\frac{1}{x-q_2^{nj}}-1\big)dx.
                                    \end{array}
\]
By evaluating $\delta_1(c_n)$ and $\delta_2(g^*(c_n))$ at $1$ in $\C_pdx$, one gets \[\delta_1(c_n)(1)=\sum_{j\geq 1}\mu(j)\frac{q_1^{jn}}{1-q_1^{jn}}dx= \sum_{j\geq 1}\sum_{k\geq 1}\mu(j)q_1^{kjn}dx=\sum_{d\geq 1}\sum_{j|d}\mu(j)q_1^{dn}dx=q_1^ndx,\] and similarly $\delta_2(g^*(c_n))(1)=q_2^ndx$.
Let $c_0$ be defined by $c_0(e_{1,j})=0$ and $c_0(e_{1,j})=1$. Then $\delta_1(c_0)=\delta_2(g^*(c_0))=\frac{dx}{x}$ and $\delta_1(c_0)(1)=\delta_2(g^*(c_0))(1)=dx$.

If $P=\sum_{n\geq 0}a_nX^n\in \mbf \Z_p[X]$, let $c_P=\sum_{n\neq 0} a_nc_n\in C(\T_1,\Z_p)$, so that
$\delta_1(c_P)(1)=P(q_1)dx$ and $\delta_2(g^*(c_P))(1)=P(q_2)dx$.

According to lemma \ref{lemordre}, $\delta_1(c_P)(1)=0$ if and only if $\delta_2(g^*(c_P))(1)=0$.
Therefore $P(q_1)=0$ if and only if $P(q_2)=0$ for every $P\in\Z_p[X]$, which implies the result.
\findem


\providecommand{\bysame}{\leavevmode\hbox to3em{\hrulefill}\thinspace}
\providecommand{\MR}{\relax\ifhmode\unskip\space\fi MR }
\providecommand{\MRhref}[2]{%
  \href{http://www.ams.org/mathscinet-getitem?mr=#1}{#2}
}
\providecommand{\href}[2]{#2}

\end{document}